\documentclass[aos,preprint]{imsart}

\usepackage{StefanTexCommandsV2}

\usepackage{amsmath, amsthm, amssymb}
\usepackage{graphicx}
\usepackage{verbatim}
\usepackage{caption}
\usepackage{subcaption}
\usepackage{fancyvrb}
\usepackage{enumerate}
\usepackage{tabularx}
\usepackage{algorithm}
\usepackage{algorithmic}
\usepackage{multirow}
\usepackage{rotating}

\usepackage{bookmark}

\RequirePackage[OT1]{fontenc}
\RequirePackage[numbers]{natbib}
\RequirePackage[colorlinks,citecolor=blue,urlcolor=blue]{hyperref}

\theoremstyle{plain}
\newtheorem{prop}{Proposition}
\newtheorem*{prop0}{Proposition}

\newtheorem{coro}[prop]{Corollary}
\newtheorem{lemm}[prop]{Lemma}
\newtheorem{theo}[prop]{Theorem}

\theoremstyle{definition}

\newtheorem{defi}{Definition}
\newtheorem{assu}{Assumption}

\floatstyle{boxed}
\theoremstyle{definition}
\newfloat{algbox}{tbp}{lop}
\newtheorem{alg}{Procedure}
\newcommand{\myalg}[3]{
\begin{center}
\vspace{-0.5\baselineskip}
\parbox{0.95\textwidth}{
\begin{alg}
\label{#1}{\textsc{ #2}}\\
\vspace{-0.5\baselineskip}
#3

\vspace{-1.5\baselineskip}
\end{alg}
}
\end{center}
}

\graphicspath{{./figures/}}

\makeatletter
\newcommand*{\maintoc}{%
  \begingroup
    \@fileswfalse
    \renewcommand*{\appendixattoc}{%
      \value{tocdepth}=-10000%
    }%
    \tableofcontents%
  \endgroup
}
\newcommand*{\appendixtoc}{
  \begingroup
    \edef\@alltocdepth{\the\value{tocdepth}}%
    \setcounter{tocdepth}{-10000}%
    \renewcommand*{\contentsname}{%
      List of Appendices}%
    \renewcommand*{\appendixattoc}{%
      \setcounter{tocdepth}{\@alltocdepth}%
    }%
    \tableofcontents%
    \setcounter{tocdepth}{\@alltocdepth}%
  \endgroup
}
\newcommand*{\appendixattoc}{}%
\g@addto@macro\appendix{%
  \clearpage%
  \phantomsection%
  \addcontentsline{toc}{section}{\appendixname}%
  \addtocontents{toc}{\protect\appendixattoc}%
}
\makeatother

\begin{document} 

\begin{frontmatter}

\title{Adaptive Concentration of Regression Trees, with Application to Random Forests}
\runtitle{Adaptive Concentration of Trees}

\runauthor{Wager and Walther}

\begin{aug}
\author{\fnms{Stefan} \snm{Wager}\corref{}\ead[label=e1]{swager@stanford.edu}}
\and
\author{\fnms{Guenther} \snm{Walther}\ead[label=e2]{gwalther@stanford.edu}}
\address{Department of Statistics \\ Stanford University \\ Stanford, CA-94305 \\ \printead{e1} \\ \printead{e2}}
\affiliation{Stanford University} 
\end{aug}

\begin{abstract}
We study the convergence of the predictive surface of regression trees and forests.
To support our analysis we introduce a notion of adaptive concentration for regression trees.
This approach breaks tree training into a model selection phase in which we pick the tree splits, followed by a model fitting phase where we find the best regression model consistent with these splits.
We then show that the fitted regression tree concentrates around the optimal predictor with the same splits:
as $d$ and $n$ get large, the discrepancy is with high probability bounded on the order of $\sqrt{\log(d)\log(n)/k}$ uniformly over the whole regression surface, where $d$ is the dimension of the feature space, $n$ is the number of training examples, and $k$ is the minimum leaf size for each tree.
We also provide rate-matching lower bounds for this adaptive concentration statement.
From a practical perspective, our result enables us to prove consistency results for adaptively grown forests in high dimensions,
and to carry out valid post-selection inference in the sense of Berk et al. [2013] for subgroups defined by tree leaves.
\end{abstract}

\end{frontmatter}

\section{Introduction}

Trees \citep{breiman1984classification} and random forests \citep{breiman2001random} are among the most widely used machine learning predictors today, with applications in a broad variety of fields such as chemistry \citep{svetnik2003random}, ecology \citep{cutler2007random,prasad2006newer}, genetics \citep{diaz2006gene,shi2004tumor}, and remote sensing \citep{ham2005investigation,pal2005random}.
While allowing for flexible predictive surfaces and complicated interactions, trees and especially random forests have proven to be surprisingly resilient to over-fitting. Unlike competing non-parametric techniques such as kernel methods or neural networks, random forests require very little tuning; experience has shown that one can often obtain good predictive models out-of-the-box with standard software like \texttt{randomForest} for \texttt{R} \citep{liaw2002classification}.
However, from the perspective of existing results, we have no particularly strong reasons to believe that forest predictions ought to be well behaved in high dimensions: The best existing convergence results for random forests either only provide fixed-dimensional asymptotic consistency guarantees \citep{breiman1984classification,scornet2014consistency}, or assume a substantially simplified training procedure where tree splits are chosen using a holdout dataset \citep{biau2012analysis,denil2014narrowing}.

\paragraph{Theoretical framework}

\begin{figure}
\includegraphics[width=0.7\textwidth]{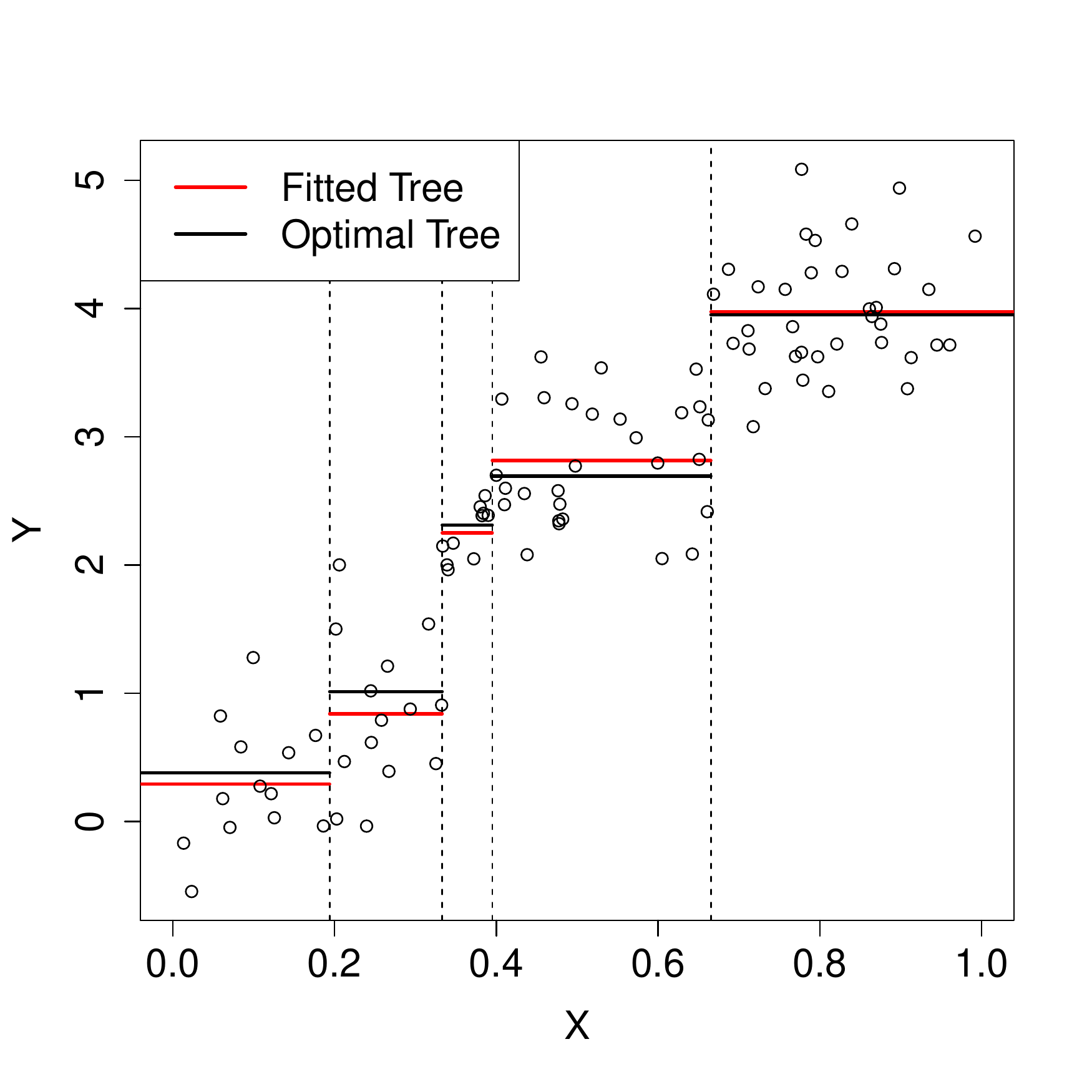}
\caption{Adaptive concentration compares the prediction surface of the fitted decision tree with that of the optimal decision tree with the same splits. Here, the splits produced by recursive partitioning are denoted by dashed vertical lines. The regression tree was fit using the \texttt{R}-package \texttt{tree} \citep{ripley2002modern}.}
\label{fig:example}
\end{figure}

The goal of this paper is to use {adaptive concentration} as a framework for describing the statistical properties of adaptively grown trees, i.e., trees that have access to the full training data while placing tree splits.
The idea of adaptive concentration is to view training trees as occurring in two stages: a model selection stage where we decide on which splits to make, and a model fitting stage where we find the best regression tree conditional on having made these splits. We then treat the splits made by the tree as fixed, and show that the fitted regression tree is not much worse than the optimal regression tree with the same splits.
In other words, we establish conditions under which sample averages over tree leaves $L$ concentrate to population averages over $L$---even if the leaves $L$ were chosen after looking at the data.
Figure \ref{fig:example} illustrates this goal for a one-dimensional tree.

Our setting is closely related to the work of \citet{berk2013valid} on valid post-selection inference, which provides convergence guarantees for estimated linear regression parameters that hold even if the regression model is selected after looking at the data.  
In order to make the connection explicit, suppose that we view trees as a form of adaptive regression, where each leaf $L$ induces an indicator feature for whether a given observation falls into the leaf. Our result then directly induces a method for constructing confidence intervals for the within-leaf mean responses that accounts for model selection in the sense of \citet{berk2013valid} (see Section \ref{sec:posi} for details).

\paragraph{Main results}

We study an asymptotic regime where the dimension $d$ of the feature space, the number $n$ of training examples and the minimum leaf size $k$ go to infinity together. Under mild regularity conditions, we show that regression trees satisfy an adaptive concentration bound that scales as $\sqrt{{\log\p{n}\log\p{d}}/{k}}$.
More specifically, with high probability and simultaneously for any leaf $L$ induced by any regression tree, the discrepancy between the sample-average response inside the leaf $L$ and the population-average response inside $L$ differ by at most $C\sqrt{{\log\p{n}\log\p{d}}/{k}}$, where $C$ is a universal constant.
In the context of Figure \ref{fig:example}, this result translates into a uniform bound on the distance between the ``fitted'' and ``optimal'' curves.
We also show that this rate of convergence is tight to within a constant factor.

Our result does not rely on major modifications to the random forest training routine, and in particular holds for variants of the CART algorithm \citep{breiman1984classification} and the original proposal of \citet{breiman2001random}. The reason there is a dependence on $n$ in the numerator of our bound is that, as the sample size grows, the trees comprising the random forest can become deeper and so the model family becomes larger.

Finally, as an application of our adaptive concentration result, we establish consistency guarantee for adaptively grown random forests---i.e., forests that do not use a holdout set for variable selection---that hold in a sparse, high-dimensional setting. To our knowledge, no directly comparable results are currently available in the literature (see Section \ref{sec:appl_consistency} for details).

\paragraph{Outline}

This paper is structured as follows. We first state and discuss our main adaptive concentration bound in Section \ref{sec:AC}, and then apply this result to provide consistency results for high-dimensional random forests in Section \ref{sec:appl_consistency}. Our theoretical work is carried out in Section \ref{sec:main_proofs}, where we develop the necessary tools to prove the results from the first sections. Section \ref{sec:lower_bounds} derives matching lower bounds for the adaptive concentration rate of regression trees.

\section{Adaptive Concentration}
\label{sec:AC}

To give adaptive concentration bounds, we first need to disambiguate the hierarchy of concepts used to build forest predictors: a {\it forest} is an ensemble of {\it trees}, each of which relies on a {\it partition} of the data generated by a splitting rule.
We begin by providing formal definitions of these quantities below; we state our main result in Section \ref{sec:results}. Throughout our analysis, we assume that we have a set of $n$ independent and identically distributed training examples $(X_i, \, Y_i)$ satisfying $X_i \in [0, \, 1]^d$ and $Y_i \in [-M, \, M]$.

\subsection{A Review of Recursive Partitioning}
\label{sec:partition}

The first concept underlying a regression tree is the splitting rule itself, which induces a partition $\Lambda$ of $[0, \, 1]^d$ into non-overlapping rectangles $L_1, \, ..., \, L_m$. We use the short-hand $L(x)$ to denote the unique element of $\Lambda$ containing $x$. We are interested in those partitions that can be obtained by {\it recursive partitioning} of the feature space \citep{breiman1984classification}.
Starting from a parent node $\nu = [0, \, 1]^d$, recursive partitioning operates by repeatedly selecting a currently unsplit node $\nu \subseteq \RR^d$, a splitting variable $j \in \{1, \, ..., \, d\}$ and a threshold $\tau \in [0, \, 1]$, and then splitting $\nu$ into two children $\nu_- = \nu \cap \{x : x_j \leq \tau\}$ and $\nu_+ = \nu \cap \{x : x_j > \tau\}$. The final leaf nodes generated by this algorithm, denoted by $L$, form a partition $\Lambda$ of $[0, \, 1]^d$.
Given our training set $\{(X_i, \, Y_i)\}$, we require the partition to be valid in the sense of Definition \ref{defi:partition}.

\begin{defi}[Valid partition]
\label{defi:partition}
A partition $\Lambda$ is {\it $\{\alpha, \, k\}$-valid} if it can by generated by a recursive partitioning scheme in which each child node contains at least a fraction $\alpha$ of the data points in its parent node for some $0 < \alpha < 0.5$, and each terminal node contains at least $k$ training examples for some $k \in \NN$. Given a dataset $\xx$, we denote the set of $\{\alpha, \, k\}$-valid partitions by $\vv_{\alpha, \, k}\p{\xx}$.
\end{defi}

The constraint that each terminal node must have at least $k$ observations is implemented by default in, e.g., \texttt{randomForest}.
Meanwhile, the requirement that each child node must incorporate at least a fraction $\alpha$ of the data in its parent---and thus that the tree cannot be excessively imbalanced---is more substantive. It is known that CART-like rules tend to split near the edges of noise features \citep{breiman1984classification,ishwaran2015effect}; thus our $\alpha$-constraint may make splits along noise features less desirable.
A similar assumption is also used by, e.g., \citep{meinshausen2006quantile}.

A partition $\Lambda$ can then be used to induce a tree predictor by averaging the responses $Y_i$ over its leaves $L_j$. In our adaptive concentration analysis, we consider two kinds of trees: valid trees \eqref{eq:leaf_mean} that are fit to the training sample, and partition-optimal trees \eqref{eq:pop_leaf_mean} that would arise if we could train a tree supported on the partition $\Lambda$ using the full population.

\begin{defi}[Valid and partition-optimal trees]
\label{defi:tree}
A valid partition induces a {\it valid tree}
\begin{equation}
\label{eq:leaf_mean}
T_\Lambda:[0, \, 1]^d \rightarrow \RR, \ \  T_\Lambda\p{x} = \frac{1}{\abs{\cb{X_i : X_i \in L(x)}}} \sum_{\cb{i : X_i \in L(x)}} Y_i.
\end{equation}
We denote the set of all $\{\alpha, \, k\}$-valid trees $T_\Lambda$ with $\Lambda \in \vv_{\alpha, \, k}\p{\xx}$ by $\tset_{\alpha, \, k}\p{\xx}$. Given a partition $\Lambda$, we also define the {\it partition-optimal tree} as
\begin{equation}
\label{eq:pop_leaf_mean}
T_\Lambda^*:[0, \, 1]^d \rightarrow \RR, \ \  T^*_\Lambda\p{x} = \EE{Y | X \in L(x)},
\end{equation}
where $(X, \, Y)$ is a new random sample from our data-generating distribution.
\end{defi}

Forests are, as their name suggests, ensembles of regression trees. Generating a regression forest involves growing multiple trees; then, the forest prediction is the average of all the tree predictions. In general, the choice of splitting splitting variables $j$ is randomized to ensure that the different trees comprising the forest are not too correlated with each other. As shown by \citet{breiman2001random}, the variance reduction of a forest in comparison with its constituent trees improves as the correlation between individual trees decreases.

\begin{defi}[Valid and partition-optimal forests]
\label{defi:forest}
For any $B \in \NN$, let $T_{\Lambda^{(1)}}, \, ..., \, T_{\Lambda^{(B)}} \in\tset_{\alpha, \, k}\p{\xx} $. Then, the average
\begin{equation}
H_{\cb{\Lambda}_1^B}: [0, \, 1]^d \rightarrow \RR, \ \ H_{\cb{\Lambda}_1^B}(x) = \frac{1}{B} \sum_{b = 1}^B T_{\Lambda^{(b)}}(x)
\end{equation}
is a {\it valid forest}; we denote the set of $\{\alpha, \, k\}$-valid forests by $\hh_{\alpha, \, k}\p{\xx}$.
The corresponding {\it partition-optimal forest} is defined as
\begin{equation}
H_{\cb{\Lambda}_1^B}^*: [0, \, 1]^d \rightarrow \RR, \ \ H_{\cb{\Lambda}_1^B}^*(x) = \frac{1}{B} \sum_{b = 1}^B T^*_{\Lambda^{(b)}}(x).
\end{equation}
When there is no risk of ambiguity, we write $H := H_{\cb{\Lambda}_1^B}$ and $H^* = H_{\cb{\Lambda}_1^B}^*$.
\end{defi}

There are many proposals for how to choose the splitting variables $j$ and the thresholds $\tau$ for trees. Our theoretical results, however, do not depend on the specific splitting rules used, and only rely on the generic structure of recursive partitioning; thus, we will not focus on specific splitting rules in this paper. For a review of how trees and forests are implemented in practice, we recommend \citet{hastie2009elements} (see Chapters 9.2 and 15).

\paragraph{Remark: Bootstrapping} One way in which our forests from Definition \ref{defi:forest} differ from the original proposal of \citet{breiman2001random} is that we do not allow the individual trees to be evaluated on bootstrap samples.\footnote{Technically, we could use a bootstrap sample to pick the partition $\Lambda$, but would then need to use the whole training set to turn $\Lambda$ into a tree $T_\Lambda$.} It seems plausible, however, that all our results should still hold even if we allow for bootstrapping, since the bootstrap is thought to have a regularizing effect on the forest and should thus reduce its ability to overfit the training data \citep{breiman1996bagging,buhlmann2002analyzing}. Studying the effect of the bootstrap on our adaptive concentration bounds and perhaps showing how it can improve adaptive concentration guarantees presents a promising avenue for further work.

\subsection{An Adaptive Concentration Bound}
\label{sec:results}

We are now ready to state our main result, given some assumptions on the problem setting.
Our first condition is a bound on the dependence of the individual coordinates $(X_i)_j$ for $j = 1, \, ..., \, d$.

\begin{assu}[Weakly dependent features]
\label{assu:unif}
We have $n$ independent and identically distributed training examples, whose features $X \in [0, \, 1]^d$ are distributed according to a density $f(\cdot)$ satisfying $\zeta^{-1} \leq f(x) \leq \zeta$ for all $x \in [0, \, 1]^d$, and some constant $\zeta \geq 1$.
\end{assu}

The above assumption is quite general. In contrast, consider the following condition:
We have features $X_i \in \RR^d$ that admit a density satisfying
\begin{equation}
\label{eq:weak_dep}
\zeta^{-1} \prod_{j = 1}^d f_j(x_j) \leq f(x) \leq \zeta \prod_{j = 1}^d f_j(x_j) \text{ for all } x \in \RR^d.
\end{equation}
Although the above may at first glance appear weaker than Assumption \ref{assu:unif}, it is actually
slightly more restrictive. 
Because trees are invariant to monotone transformations of the features $(X_i)_j$, we can without
loss of generality rescale the features such that $X_i \in [0, \, 1]^d$ with uniform marginals. Applying
this transformation to \eqref{eq:weak_dep} would yield Assumption \ref{assu:unif}, along with a uniformity
constraint on the marginals of $f(\cdot)$.

Next, although we allow for the minimum leaf size $k$ to be quite small, we still need for it to grow with $n$. Note that since $k \leq n$, the assumption below implicitly requires that $\log(d) \ll n/\log(n)$.

\begin{assu}[Minimum leaf size]
\label{assu:scaling}
The minimum leaf-size $k$ grows with $n$ at a rate bounded from below by
\begin{equation}
\label{eq:k_lb}
\limn \frac{\log\p{n} \, \max\cb{\log(d), \, \log\log(n)}}{k} = 0.
\end{equation}
\end{assu}

The following theorem is our main result on the adaptive concentration of decision trees.
This result implies that, in practical data analysis, we can treat the fitted prediction function $T_\Lambda$ as a good approximation to the optimal tree $T_\Lambda^*$ supported on the partition $\Lambda$.
This result requires no continuity assumptions on the conditional mean function $\EE{Y \cond X = x}$.

\begin{theo}
\label{theo:main}
Suppose that we have $n$ training examples $(X_i, \, Y_i) \in [0, \, 1]^d \times [-M, \, M]$ satisfying Assumption \ref{assu:unif}, and that we have a sequence of problems with parameters $(n, \, d, \, k)$ satisfying Assumption \ref{assu:scaling}.
Then, sample averages over all possible valid partitions concentrate around their expectations with high probability:
\begin{align}
\label{eq:main_complete}
\begin{split}
&\lim_{n, \, d, \, k \rightarrow \infty} \mathbb{P}\vast[\sup_{x \in [0, \, 1]^d, \, \Lambda \in \vv_{\alpha, \, k}} \left|T_\Lambda\p{x} - T^*_\Lambda\p{x} \right|
\\  &\ \ \ \ \ \ \ \ \ \
\leq 9 M  \sqrt{\frac{ \log \p{n/k} \p{\log(dk) + 3\log\log(n)}}{\log\p{\p{1 - \alpha}^{-1}}}} \, \frac{1}{\sqrt{k}} \vast] = 1.
\end{split}
\end{align}
In a moderately high-dimensional regime with $\liminf d/n > 0$, the above bound simplifies to
\begin{align}
\label{eq:main}
\mathbb{P}\sqb{\sup_{x \in [0, \, 1]^d, \, \Lambda \in \vv_{\alpha, \, k}} \left|T_\Lambda\p{x} - T^*_\Lambda\p{x} \right|
\leq 9 M  \sqrt{\frac{ \log \p{n} \log \p{d}}{\log\p{\p{1 - \alpha}^{-1}}}} \, \frac{1}{\sqrt{k}}} \rightarrow 1.
\end{align}
\end{theo}

This result appears to be remarkably strong. As a baseline, suppose we just selected a single tree $T \in \tset_{\alpha, \, k}\p{\xx}$ non-adaptively, i.e., without looking at the labels $Y_i$. Then, a simple Hoeffding bound where we take $n$ as a crude upper for the total number of leaves shows that
\begin{equation}
\label{eq:non_adaptive}
\PP{\sup_{x \in [0, \, 1]^d} \left|T\p{x} - T^*\p{x} \right| \leq M \sqrt{\frac{2.1\log\p{n}}{k}}} \rightarrow 1.
\end{equation}
Moreover, assuming that $k$ is not too large, say $k \leq \sqrt{n}$, and that the tree is fully grown to depth $k$, the bound \eqref{eq:non_adaptive} is essentially tight.

Comparing \eqref{eq:main} with \eqref{eq:non_adaptive}, we see that a uniform concentration bound valid for {all} possible trees in $\tset_{\alpha, \, k}\p{\xx}$ is only a factor $\oo(\sqrt{\log(d)})$ weaker than the best concentration bound we could hope for with a {single} tree. In other words, the ``cost'' of adaptively searching over all valid trees in high dimensions is surprisingly low, and only scales with $\sqrt{\log(d)}$.

In Section \ref{sec:lower_bounds}, we provide a rate-matching lower bound to Theorem \ref{theo:main} that holds for $d = n^r$ for some $r \geq 1$. Thus, it is not possible to substantially improve on the above result without restricting the class of considered trees.

Finally, because Theorem \ref{theo:main} holds simultaneously for all trees $T \in \tset$, we note that \eqref{eq:main_complete} and \eqref{eq:main} also induce adaptive concentration bounds for random forests. The following result is conceptually related to the work of \citet{biau2008consistency}, who established conditions under which ensembles inherit consistency properties of their base classifiers, despite potentially allowing for greater representational power.

\begin{coro}
\label{coro:main_rf}
Under the conditions of Theorem \ref{theo:main}, and assuming that $\liminf d/n > 0$, we find that 
\begin{equation}
\label{eq:main_rf}
\mathbb{P}\sqb{\sup_{x \in [0, \, 1]^d} \left|H\p{x} - H^*\p{x} \right|
\leq 9 M  \sqrt{\frac{ \log \p{n} \log \p{d}}{\log\p{\p{1 - \alpha}^{-1}}}} \, \frac{1}{\sqrt{k}}} \rightarrow 1.
\end{equation}
simultaneously for all valid forests $H$ constructed as in Definition \ref{defi:forest}.
\end{coro}

Several authors have studied the theoretical properties of trees and forests in order to explain their performance
\citep{arlot2014analysis,biau2012analysis,biau2008consistency,breiman2004consistency,denil2014narrowing,gordon1984almost,
lin2006random,meinshausen2006quantile,olshen2007tree,scornet2014asymptotics,
scornet2014consistency,scott2006minimax,wager2014asymptotic}.
In particular, \citet{scornet2014consistency} use low-dimensional asymptotics prove that Breiman's original forests are consistent in terms of their predictive risk, while \citet{biau2012analysis} and \citet{denil2014narrowing} discuss the properties of some random forests where trees are grown using only a holdout set---i.e., without looking at the training data used for predictions.
To our knowledge, however, our above result is the first convergence guarantee for the predictive surface of tree-based regression that holds in an asymptotic regime where $n$ and $d$ go to infinity together, and does not require the use of a holdout set.

\subsection{Adaptive Concentration as Post-Selection Inference}
\label{sec:posi}

We can also interpret Theorem \ref{theo:main} directly as a result about statistical inference for adaptive regression. To motivate this view, suppose that we build a tree-based model to do subgroup analysis on data from a clinical trial with the intent of, e.g., prioritizing treatment to patients in leaves with the largest estimated treatment effect \citep{athey2015machine,su2009subgroup}. In this case, it is crucial to have valid {post-selection} confidence intervals for the treatment effect within each leaf that account for the adaptive model search done by the tree.\footnote{Given an unconfoundedness assumption \citep{rosenbaum1983central}, results about regression trees and forests can directly be adapted to the problem of heterogeneous treatment effect estimation; see Section 4 of \citet{wager2014asymptotic} for an example.}

At a high level, the problem of post-selection inference for trees is a special case of post-selection inference for adaptive least-squares regression. In the classical adaptive regression setup, we start with a design $X \in \RR^{n \times d}$ and a response $Y \in \RR^n$; the statistician then selects a model $\mm \subseteq \cb{1, \, ..., \, d}$ and provides estimates of the form $\hy = A \hbeta_A$, where $A := X_\mm \in \RR^{n \times \abs{\mm}}$ is the matrix comprised of the columns of $X$ contained in the set $\mm$ \citep[e.g.,][]{berk2013valid,fithian2014optimal,lee2013exact}.

In contrast, as in Definition \ref{defi:tree}, a valid tree first generates a partition $\Lambda = \{L_j\}_{j = 1}^m$, and then averages the responses $Y_i$ inside each leaf $L_j$. Formally, this procedure is equivalent to creating a design matrix
$$ A \in \cb{0, \, 1}^{n \times m}, \text{ where } A_{ij} = 1\p{\cb{X_i \in L_j}}, $$
and then running linear regression on $A$ and the response vector $Y$. In other words, regression trees are a form of adaptive regression where instead of subsetting a list of available variables, we build a design matrix by considering indicator functions over leaves.

For any such adaptively-constructed design $A$, we can define the optimal regression vector
$$ \beta^*_A = \p{A^\top A}^{-1} A^\top \mu_i \text{ where } \mu_i = \EE{Y_i \cond A_i}; $$
it is then natural to ask about the discrepancy between $\hbeta_A - \beta^*_A$,
where $\hbeta_A = (A^\top A)^{-1} A^\top Y$ is the ordinary least squares estimate in the adaptive regression model.
In the case of subgroup analysis with trees, $\hbeta_A$ corresponds to the measured mean responses in each tree leaf, whereas $\beta^*_A$ encodes the population averages over the same leaves.

The problem of bounding the gap between $\hbeta_A - \beta^*_A$ for sparse linear regression, i.e., where $A = X_\mm$ for some feature subset $\mm$, has been studied in detail by \citet{berk2013valid}. The authors show that it is possible to honestly account for arbitrary model selection by inflating standard confidence intervals for linear regression by a {PoSI constant} $K$ that depends on $X$. For orthogonal designs, the optimal PoSI constant scales as $K \sim \sqrt{2\log d}$; however, there also exist designs for which $K \geq \sqrt{d}/2$. Moreover, computing the PoSI constant is in general difficult; the method used by \citet{berk2013valid} only works for $d \leq 20$.

In contrast, our Theorem \ref{theo:main} shows that, for adaptively grown regression trees,
$\lVert\hbeta_A - \beta^*_A\rVert_{\infty} = \oo_p(\sqrt{\log(n)\log(d)/k})$;
in other words, there exists a universal PoSI constant for trees on the order of $\sqrt{\log(n)\log(d)}$.
This is somewhat worse than the PoSI correction for linear regression with an orthogonal design, but is exponentially better than general-case performance for sparse linear regression.

\subsection{Adaptive Concentration vs. Generalization Bounds}

There is a long tradition of studying the convergence of regression trees using arguments in the style of \citet{vapnik1971uniform}---including the work of \citet{breiman1984classification} and, more recently,
\citep{bartlett2003rademacher,devroye1996probabilistic,lugosi1996consistency,mansour2000generalization,nobel1996histogram}. 
The goal of these papers is to control the generalization error of regression trees, i.e., the difference between the training error and the expected test error of the tree.

In our context, we can use \citep{mansour2000generalization} to obtain the following generalization bound:
\begin{align}
\label{eq:tree_genbound}
\begin{split}
\sup_{T \in \tset_{\alpha, \, k}} \ &\frac{1}{n}\sum_{i = 1}^{n} \p{Y_i - T(X_i)}^2 - \EE{\p{Y - T(X)}^2} \\
&= \oo_P\p{\sqrt{\frac{\#\cb{\text{leaves}}  \p{\log\p{d} + \log\p{n}}}{n}}} \\
&= \oo_P\p{\sqrt{\frac{\log\p{d} + \log\p{n}}{k}}}.
\end{split}
\end{align}
Strikingly, the rate of convergence in \eqref{eq:tree_genbound} is better than the adaptive concentration rate obtained in Theorem \ref{theo:main}. Now, given our matching lower bound (see Theorem \ref{theo:lowerbound}), we know that our rates are optimal in moderately high dimensions; thus, good adaptive concentration must be fundamentally more difficult than good generalization.

The reason for this discrepancy is that generalization bounds only seek to control the global performance of regression trees, whereas adaptive concentration requires local control of the tree regression surface. 
To give a concrete example, suppose that $d = n^r$ for some $r > 0$. Then, the bound \eqref{eq:tree_genbound} implies that the empirical risk of a regression tree will be consistent for the test set error of the same tree provided that $k \gg \log(n)$. Meanwhile, our adaptive concentration analysis implies that there exist trees with bad leaves unless $k \gg \log(n)^2$; however, these bad leaves will be rare enough as to not affect the overall test set error of the tree.

Generalization bounds vs. adaptive concentration bounds can both be relevant in different practical applications. If all we care about is the test error of a tree, then the bound \eqref{eq:tree_genbound} is more useful than Theorem \ref{theo:main}. Conversely, if we want to use a regression tree to identify outlying sub-populations, then we need guarantees as in Theorem \ref{theo:main}, whereas the rate \eqref{eq:tree_genbound} might be misleading.

\section{Consistency of High-Dimensional Random Forests}
\label{sec:appl_consistency}

As one practical application of Theorem \ref{theo:main}, we obtain consistency guarantees for random forests in a high-dimensional setting. We work in a regime where $n, \, d \rightarrow \infty$, but the conditional mean function $\EE{Y \cond X = x}$ only depends on a small number of covariates:

\begin{assu}[Sparse signal]
\label{assu:sparse}
There is a signal set $\qq \in \cb{1, \, ..., \, d}$ of size $\abs{\qq} \leq q$ such that the set of random variables $\{(X_i)_j : j \not\in \qq\}$ is jointly independent of $Y_i$ and the set $\{(X_i)_j : j \in \qq\}$.
\end{assu}

We study a simple CART-like regression forest that can consistently estimate sparse signals in high dimensions, described in Procedure \ref{alg:guessandcheck}. Effectively, the algorithm is an extreme form of the procedure of \citet{breiman2001random}, where each splitting variable is chosen uniformly at random from $\{1, \, ..., \, d\}$. However, in a break from classical regression trees, our algorithm uses Theorem \ref{theo:main} to test the significance of a candidate split: if we have never split on a variable $j$ yet, then our tree will only split along $j$ if it leads to a large enough improvement in mean-squared error. Once the tree has split along $j$ once, then this feature is ``unlocked'' and can be subsequently used without testing.\footnote{Running a hypothesis test before accepting a split for a regression tree has also been considered by \citet{zeileis2008model}; however, the paper does not provide any formal consistency guarantees.}
At a high level, our construction relies on a guarantee from Theorem \ref{theo:main} that no noise feature $j$ will ever appear significant enough to get unlocked at any stage of the forest-generation process.

\begin{algbox}[t]
\myalg{alg:guessandcheck}{Guess-and-Check Forest}{

Input: $n$ training examples of the form $(X_i, \, Y_i)$, a minimum leaf size $k$, and a balance parameter $\alpha$.

\vspace{0.5\baselineskip}

Guess-and-check trees recursively apply the following splitting procedure until no more splits are possible, i.e., until all terminal nodes contain less than $2k$ training examples or there were no possible splits satisfying \eqref{eq:split_thresh} to begin with.

\begin{enumerate}
\item Select a currently unsplit node $\nu$ containing at least $2k$ training examples.
\item Pick a candidate splitting variable $j \in \cb{1, \, ..., \, d}$ uniformly at random.
\item Pick the minimum squared error splitting point \smash{$\htheta$}. More specifically,
\begin{align*}
\htheta = \ \ & {\bf argmax} \ \ell\p{\theta} := \frac{4 \, N^-\p{\theta} \, N^+\p{\theta}}{\p{N^-\p{\theta} + N^+\p{\theta}}^2} \, \Delta^2\p{\theta} \\
& \text{\bf such that }  \theta = (X_i)_j \text{ for some } X_i \in \nu\\
& \ \ \ \ \ \ \ \ \alpha \abs{\cb{i : X_i \in \nu}}, \, k \leq N^-\p{\theta}, \, N^+\p{\theta} \\
&\text{\bf where } \Delta\p{\theta} = \hspace{-6mm} \sum_{\cb{i : X_i \in \nu(x), \, (X_i)_j > \theta}} \hspace{-6mm}  Y_i \, \big/ \, N^+ \hspace{2mm} -  \hspace{-6mm} \sum_{\cb{i : X_i \in \nu(x), \, (X_i)_j \leq \theta}} \hspace{-6mm}  Y_i \, \big/ \, N^-, \\
& \ \ \ \ \ \ \ \ \ \ N^-\p{\theta} = \abs{\cb{i : X_i \in \nu, \, (X_i)_j \leq \theta}}, \\
& \ \ \ \ \ \ \ \ \ \ N^+\p{\theta} = \abs{\cb{i : X_i \in \nu, \, (X_i)_j > \theta}}.
\end{align*}
\item If {\bf either} there has already been a successful split along variable $j$ for some other node {\bf or}
\begin{equation}
\label{eq:split_thresh}
\ell\p{\htheta} \geq \p{2 \times 9 M \, \sqrt{\frac{\log\p{n} \log\p{d}}{k \, \log\p{\p{1 - \alpha}^{-1}}}}}^2,
\end{equation}
the split succeeds and we cut the node $\nu$ at $\htheta$ along the $j$-th variable; if not, we do not split the node $\nu$ this time.
\end{enumerate}

A guess-and-check forest is the average of $B$ independently generated guess-and-check trees.

\vspace{0.5\baselineskip}

}
\end{algbox}

To guarantee consistency of guess-and-check forests, we still need two additional conditions to hold. Assumption \ref{assu:monotone} is a weaker version of an assumption made by \citet{scornet2014consistency} to ensure that trees can eventually accumulate evidence enabling them to split along signal coordinates. The continuity assumption is standard for consistency proofs \citep[e.g.,][]{biau2012analysis,biau2008consistency,meinshausen2006quantile,scornet2014consistency}.

\begin{assu}[Monotone signal]
\label{assu:monotone}
There is a minimum effect size $\beta > 0$ and a set of sign variables $\sigma_j \in \cb{\pm 1}$ such that, for all $j \in \qq$ and all $x \in [0, \, 1]^d$,
\begin{align*}
&\sigma_j \, \bigg( \EE{Y_i \cond (X_i)_{-j} = x_{-j}, (X_i)_j > \frac{1}{2}} \\
&\ \ \ \ \ \ \ \ -  \EE{Y_i \cond (X_i)_{-j} = x_{-j}, (X_i)_j\leq \frac{1}{2}} \bigg) \geq \beta,
\end{align*}
where $x_{(-j)} \in [0, \, 1]^{d - 1}$ denotes the vector containing all but the $j$-th coordinate of $x$.
\end{assu}

\begin{assu}[Continuity]
\label{assu:lip}
The function $\EE{Y \cond X = x}$ is Lipschitz-continuous in $x$.
\end{assu}

Given this setup, our first result is a uniform consistency result for guess-and-check forests in high dimensions; recall that, as in Theorem \ref{theo:main}, we assume that the minimum leaf size $k$ grows faster than $\log(n)\log(d)$. We take the signal dimension $q$ to be fixed.

\begin{theo}
\label{theo:unif_consistent}
Under the conditions of Theorem \ref{theo:main} with $\liminf d/n > 0$, suppose that $\hy(x)$ are estimates for $\EE{Y \cond X = x}$ obtained using a guess-and-check forest (Procedure \ref{alg:guessandcheck}). Suppose, moreover, that Assumptions \ref{assu:sparse}, \ref{assu:monotone}, and \ref{assu:lip} hold. Then,
$$ \lim_{n \, d, \, k \rightarrow \infty} \sup_{x \in [0, \, 1]^d} \abs{\hy\p{x} - \EE{Y \cond X = x}} = 0. $$
\end{theo}

To our knowledge, this is the first pointwise consistency result for adaptively grown random forests in high dimensions. \citet{scornet2014consistency} prove that random forests are consistent; however, their proof only works in a regime where $d$ is fixed while $n \rightarrow \infty$.
Meanwhile, in high dimensions, \citet{biau2012analysis} establishes a rate of convergence for random forests that only depends on the signal dimension $q$ and not on $d$. However, his proof assumes a form of data splitting, where we first effectively do variable selection on a holdout sample, and then grow non-adaptive trees on the training sample without looking at the responses $Y$.

We also consider the predictive error of our guess-and-check forests. To do so, we focus on a special variant of guess-and-check trees with ``$\alpha = 0.5$'', meaning that $\htheta$ in step 3 of Procedure \ref{alg:guessandcheck} is always set to the median of the $(X_i)_j$ with $X_i \in \nu$. The idea of building trees by splitting along leaf medians goes back to at least \citet{devroye1996probabilistic}, and the behavior of median forests in low dimensions has been studied in detail by \citet{duroux2016impact}.

The following result hinges on showing that our median guess-and-check forest, with examples in $[0, \, 1]^d$, converges as fast as a standard median forest with examples in $[0, \, 1]^q$. We can thus recover the rates of convergence of \citet{duroux2016impact} that only depend on $q$. In order to apply the result of \citet{duroux2016impact}, we need to require the $X_i$ to be uniformly distributed over $[0, \, 1]^d$, i.e., $\zeta = 1$ in Assumption \ref{assu:unif}.

\begin{theo}
\label{theo:l2_consistent}
Set $\alpha = 0.5$, and define $\xi = 1/(1 - 3/(4q))$. Suppose that the conditions of Theorem \ref{theo:main} with $\liminf d/n > 0$ and Assumptions \ref{assu:sparse}, \ref{assu:monotone} and \ref{assu:lip} hold, and moreover that $X_i$ is uniformly distributed over $[0, \, 1]^d$. Then, the excess error rate of the guess-and-check tree (Procedure \ref{alg:guessandcheck}) satisfies:
\begin{equation}
\label{eq:biau_rate}
\EE{\p{\hy\p{X} - \EE{Y \cond X}}^2} = \oo\p{n^{-\frac{\log\p{\xi}}{\log\p{2 \, \xi}}}},  \ \text{given} \  k \asymp n^\frac{\log\p{\xi}}{\log\p{2 \,\xi}}.
\end{equation}
\end{theo}

The result below is a direct analogue to the result of \citet{biau2012analysis}, except that our trees are adaptive, i.e., grown as usual by considering the $Y$ sample; whereas his result required a holdout sample for variable selection.
As noted by Biau, as soon as $q < 0.54 \, d$, the rate \eqref{eq:biau_rate} is better than the standard rate for non-parametric estimation of Lipschitz functions, i.e., $n^{-2/(d+2)}$.

Finally, we note that Assumption \ref{assu:monotone} was only used as a crude tool to ensure that the guess-and-check forest succeeds at splitting on all the signal variables; given Assumptions \ref{assu:sparse} and \ref{assu:lip}, our proof shows that guess-and-check forests will be consistent as long as the trees in fact split on the signal variables. Empirically, greedily trained regression trees have been found to be powerful under substantially weaker conditions than Assumption \ref{assu:monotone}, so we expect a relaxation to be possible; however, we leave this to further work.

\subsection{Proof Sketch}
\label{sec:consistency}

In this section, we briefly outline the main steps in proving our forest consistency results using Theorem \ref{theo:main}.
We begin by showing that, as $n, \, d \rightarrow \infty$, guess-and-check trees (Procedure \ref{alg:guessandcheck}) never split on noise variables. We emphasize that the following result holds simultaneously for {all} random realizations of the guess-and-check forest.

\begin{lemm}
\label{lemm:no_bad_split}
Under the conditions of Theorem \ref{theo:main} with $\liminf d/n > 0$, suppose moreover that Assumption \ref{assu:sparse} holds. Let $\pi_{bad}$ be the probability that any guess-and-check tree ever splits on a noise variable $j \not\in \qq$. Then,
$\pi_{bad} = \oo(1/\sqrt{n})$.
\end{lemm}

In order to achieve consistency, we also need to make enough splits along the signal variables. The following lemma guarantees that, given Assumption \ref{assu:monotone}, a guess-and-check tree will almost always split on a signal variable when it gets a chance to do so.

\begin{lemm}
\label{lemm:good_split}
Suppose Assumptions \ref{assu:sparse} and \ref{assu:monotone} hold. For any variable $j$ for any $j \in \qq$, let $\pi_j$ be the probability that the first time any guess-and-check tree tries to split along $j$, the split succeeds. Then
$\pi_j = 1 - \oo(1/\sqrt{n})$.
\end{lemm}

Given these two lemmas, we see that a $d$-dimensional guess-and-check tree with a $q$-dimensional signal effectively behaves like a $q$-dimensional tree. Then, to prove Theorems \ref{theo:unif_consistent} and \ref{theo:l2_consistent}, it suffices to use existing results about the consistency of random forests in low dimensions; we in particular build on the work of \citet{duroux2016impact} and \citet{meinshausen2006quantile}.

\section{Theoretical Development}
\label{sec:main_proofs}

We now develop the technical machinery required to prove Theorem \ref{theo:main}.
The bulk of our work goes into bounding large deviations of the process
\begin{equation}
\label{eq:emp_proc}
 \frac{1}{\abs{\cb{i : X_i \in L}}} \sum_{\cb{i : X_i \in L}} Y_i - \EE{Y \cond X \in L},
\end{equation}
where $L$ ranges over the set $\law_{\alpha,\,k}$ of all possible leaves of a valid partition $\Lambda \in \vv_{\alpha,\,k}$.
Our argument proceeds as follows.

First, in Section \ref{sec:rectangles}, we construct a parsimonious set of rectangles $\mathcal{R}$ such that any leaf $L \in \law_{\alpha,\,k}$ can be well-approximated by a rectangle $R \in \mathcal{R}$ under Lebesgue measure $\lambda(\cdot)$. More specifically, for any large enough $L \in \law_{\alpha,\,k}$, we require that there exist $R_-, \, R_+ \in \mathcal{R}$ such that $R_- \subseteq L \subseteq R_+$, and $e^{-1/\sqrt{k}} \, \lambda(R_+) \leq \lambda(L) \leq e^{1/\sqrt{k}} \, \lambda(R_-)$.

Then, in Section \ref{sec:concentration}, we establish concentration bounds for the process \eqref{eq:emp_proc} that only depend on the cardinality of the approximating set $\mathcal{R}$. The main step here is in showing that the rectangles $R_-$ and $R_+$ constructed above are also good approximations to $L$ under the {empirical} measure; for example, if $R_-$ is an inner approximation to $L$, then $R_-$ cannot contain too many fewer training examples than $L$.
The technical results from Sections \ref{sec:rectangles} and \ref{sec:concentration} directly yield Theorem \ref{theo:main}, as explained in Section \ref{sec:main_result_proof}.
Proofs are given in the appendix.

\subsection{Notation}
\label{sec:notation}

Throughout our analysis, we assume that we have $n$ labeled independent and identically distributed training examples $(X_i, \, Y_i) \in [0, \, 1]^d \times \RR$. We denote rectangles $R \in [0, \, 1]^d$ by
\begin{equation}
R = \bigotimes_{j = 1}^d \sqb{r_j^-, \, r_j^+}, \where 0 \leq r_j^- < r_j^+ \leq 1 \text{ for all } j = 1, \, ..., \, d,
\end{equation}
writing the Lebesgue measure of $R$ as $\lambda\p{R}$
\begin{equation}
\label{eq:lebesgue}
\lambda\p{R} = \prod_{j = 1}^d \p{r_j^+ - r_j^-}.
\end{equation}
Recall that, by Assumption \ref{assu:unif}, the features $X_i$ have a distribution $f(\cdot)$ satisfying $\zeta^{-1} \leq f(x) \leq \zeta$ for all $x \in [0, \, 1]^d$. Given this setting, we write the expected fraction of training examples falling inside $R$ as $\mu(R)$,  and the number of training examples $X_i$ inside $R$ as $\#R$:
\begin{equation}
\label{eq:counts}
\mu\p{R} =\int_R f(x) \ dx, \ \  \#R = \abs{\cb{i : X_i \in R}}.
\end{equation}
Notice that, marginally, $\#R \sim \text{Binomial}\p{n,\, \mu\p{R}}$. For any rectangle $R$, we define its support as
\begin{equation}
\label{eq:support}
 S(R) = \cb{j \in 1, \, ..., \, d : r_j^{-} \neq 0 \text{ or } r_j^+ \neq 1};
\end{equation}
these are the features used in defining $R$. Finally, we write $\law_{\alpha,\,k}$ for the set of all possible leaves associated with a valid partition $\Lambda \in \vv_{\alpha,\,k}$.  This notation is summarized in Table \ref{tab:notation}.

\begin{table}[t]
\caption{Summary of notation.}
\label{tab:notation}
\begin{tabular}{|c|c|}
\hline
$n$ & Number of training examples \\
$d$ & Dimension of feature space: $X \in [0, \, 1]^d$ \\
$k$ & Minimum leaf size \\
\hline
$\alpha$ & Bound on the allowable imbalance in recursive partitioning \\
$\vv_{\alpha, \, k}$ & Set of $\cb{\alpha,\, k}$-valid partitions; see Definition \ref{defi:partition} \\
$\law_{\alpha, \, k}$ & Set of all leaves generated by $\cb{\alpha,\, k}$-valid partitions \\
$\tset_{\alpha, \, k}$ & Set of $\cb{\alpha,\, k}$-valid trees; see Definition \ref{defi:tree} \\
$\hh_{\alpha, \, k}$ & Set of $\cb{\alpha,\, k}$-valid forests; see Definition \ref{defi:forest} \\
\hline
$\zeta$ & Bound on the dependence structure of $X$; see Assumption \ref{assu:unif} \\
$\lambda(R)$ & Lebesgue measure of a rectangle $R$; see \eqref{eq:lebesgue} \\
$\mu(R)$ & Expected fraction of training examples in a rectangle $R$; see \eqref{eq:counts} \\
$\#R$ & Number of training examples inside a rectangle $R$; see \eqref{eq:counts} \\
$S(R)$ & Support of a rectangle $R$; see \eqref{eq:support} \\
\hline
\end{tabular}
\end{table}

\subsection{A Set of Approximating Rectangles}
\label{sec:rectangles}

Our first result effectively bounds the complexity of the space of rectangles over the unit cube by showing how all such rectangles can be well-approximated using an economical set of rectangles $\rr$. We detail a constructive characterization of $\rr$ in Section \ref{sec:construction}; this construction is a generalization of the one used by \citet{walther2010optimal} to study the asymptotics of the multidimensional scan statistic \citep{kulldorff1997spatial}.

\begin{theo}
\label{theo:approx_set}
Let $S \in \cb{1, \, ..., \, d}$ be a set of size $|S| = s$, and let $w, \, \varepsilon \in (0, \, 1)$. Then, there exists a set of rectangles $\rr_{S, \, w, \, \varepsilon}$ such that the following properties hold.
Any rectangle $R$ with support $S(R) \subseteq S$ and of volume $\lambda\p{R} \geq w$ can be well approximated by elements in $\rr_{S, \, w, \, \varepsilon}$ from both above and below in terms of Lebesgue measure. Specifically, there exist rectangles $R_-, \, R_+ \in \rr_{S, \, w, \, \varepsilon}$ such that
\begin{align}
\label{eq:eps_app}
&R_- \subseteq R \subseteq R_+, \ \eqand \
e^{-\varepsilon} \lambda\p{R_+} \leq \lambda\p{R} \leq e^\varepsilon \lambda\p{R_-}.
\end{align}
Moreover, the set $\rr_{S, \, w, \, \varepsilon}$ has cardinality bounded by
\begin{equation}
\abs{\rr_{S, \, w, \, \varepsilon}} \leq \frac{1}{w} \p{\frac{8 s^2}{\varepsilon^2} \p{1 + \log_2 \left\lfloor \frac{1}{w} \right \rfloor}}^s  \cdot \p{1 + \oo\p{\varepsilon}}.
\end{equation}
\end{theo}

\newcommand{\rrset}{\rr_{s, \, w, \, \varepsilon}}

In order to approximate all possible $s$-sparse rectangles, we use the set
\begin{equation}
\label{eq:union}
\rrset = \cup_{|S| = s} \rr_{S, \, w, \, \varepsilon}
\end{equation}
of size
\begin{equation}
\abs{\rrset} \leq \binom{d}{s} \, \frac{1}{w} \p{\frac{8 s^2}{\varepsilon^2} \p{1 + \log_2 \left\lfloor \frac{1}{w} \right \rfloor}}^s  \cdot \p{1 + \oo\p{\varepsilon}}.
\end{equation}
Now, given our tree construction as encoded in Definition \ref{defi:partition},
each child node must be smaller than its parent by at least a factor $1 - \alpha$;
thus, for any $L \in \law$,
$ \#L \leq \p{1 - \alpha}^{\abs{S(L)}} n $,
and so
\begin{align}
\abs{S(L)}
\leq \frac{\log\p{n / \#L}}{\log\p{1 / \p{1 - \alpha}}} 
\leq  \frac{\log\p{n / k}}{\log\p{1 / \p{1 - \alpha}}}.
\end{align}
This implies that, by setting
\begin{equation}
\label{eq:s_def}
s = \left\lfloor \frac{\log\p{n/k}}{\log\p{1/\p{1 - \alpha}}} \right\rfloor + 1,
\end{equation}
we can use $\rrset$ to approximate all possible tree leaves $L \in \law_{\alpha, \, k}$
to within error $\varepsilon$ under Lebesgue measure $\lambda(\cdot)$,
provided the leaves have volume $\lambda(L) \geq w$.
We end this section with a useful bound on the size of the approximating set $\rrset$. 

\begin{coro}
\label{coro:set_size_bound}
Suppose that we set
\begin{align}
\label{eq:param_choices}
w = \frac{1}{2 \zeta}\, \frac{k}{n}, \ \ \varepsilon = \frac{1}{\sqrt{k}}, \, \eqand \, s = \left\lfloor\frac{\log\p{n/k}}{\log\p{\p{1 - \alpha}^{-1}}}\right\rfloor + 1,
\end{align}
where $0 < \alpha < 0.5$ and $\zeta \geq 1$ are fixed constants. Then,
\begin{align}
\label{eq:set_size_bound}
\begin{split}
&\log\p{\abs{\rr_{s, \, w, \, \varepsilon}}} \leq  \frac{\log\p{n/k} \p{\log\p{d k} + 3\log\log(n)}}{\log\p{\p{1 - \alpha}^{-1}}} 
\\  &\ \ \ \ \ \ \ \ \ \
 +  \oo\p{\log\p{\max\cb{n, \, d}}}.
 \end{split}
\end{align}
\end{coro}

\subsubsection{Constructing Approximating Rectangles}
\label{sec:construction}

Without loss of generality, we can take $S = \{1, \, ...,\, s\}$; thus, our job is to $\varepsilon$-approximate all rectangles $R \in [0, \, 1]^s$ of volume at least $w$.
When $s = 1$, it is easy to verify that we can construct an approximating set containing on the order of $w^{-2}$ elements that $\varepsilon$-approximate all possible intervals of length greater than $w$: we can build such a set by, e.g., considering all rectangles of the form $[a \cdot w \varepsilon/2, \, b \cdot w \varepsilon/2]$ where $a$ and $b$ are integers.

A naive extrapolation of this idea may suggest that, as $s$ grows, the number of required rectangles scales as $w^{-2s}$: this is what we would get by varying all the parameters $r_j^-$ and $r_j^+$ freely. However, as shown by the construction below, this guess is much too pessimistic. The reason for this is that the volume constraint $\lambda(R) = \prod_{j = 1}^s (r_j^+ - r_j^-) \geq w$  becomes more and more stringent as the dimension $s$ grows, because every dimension along which $r_j^- \not\approx 0$ or $r_j^+ \not\approx 1$ geometrically cuts the size of $\lambda(R)$. For example, we can immediately verify that if $\lambda(R) \geq w$, then $r_j^+ - r_j^- \leq 0.5$ can can hold for at most $\log_2 (w^{-1})$ coordinates.

The construction below exploits the intuition that at most a few coordinates can be active on a small scale. Generalizing ideas from \citep{walther2010optimal}, we define $\rr$ as the set of all rectangles of the form $R = \bigotimes_{j = 1}^s  [r_j^-, \, r_j^+ ]$, with
\begin{align}
&r_j^- = a_j 2^{\tau_j - 1} \, \frac{w \varepsilon}{s} \,\eqand\,
r_j^+ = \min\cb{1, \, r_j^- + w 2^{\tau_j} + b_j 2^{\tau_j - 1} \, \frac{w \varepsilon}{s}},
\end{align}
such that
\begin{align}
\label{eq:counter}
&a_j \in 0, \, 1, \, .., \, \left\lfloor 2^{1 - \tau_j} \frac{s}{w \varepsilon} \right\rfloor, \
b_j \, \in \, 0, \, 1, \, .., \, \left\lceil \frac{2s}{\varepsilon} \right\rceil, \\
\label{eq:counter2}
&\tau_j \in 0, \, 1, \, ..., \, \lfloor \log_2 w^{-1} \rfloor, \, \eqand\,
\sum_{j = 1}^s \tau_j\geq \p{s - 1}\log_2\p{\frac{1}{w}} - s.
\end{align}
In this construction, the $j$-th interval $[r_j^-, \, r_j^+]$ is on the scale $w \, 2^{\tau_j}$. The observation that only a few coordinates $j$ can be active on a small scale is encoded in the lower bound \eqref{eq:counter2} on $\sum_{j = 1}^s \tau_j$. The lemma below confirms that this this approximating set is valid.

\begin{lemm}
\label{lemm:eps_approx}
Given any rectangle $R$ with support $S(R) \subseteq S$ and volume $\lambda\p{R} \geq w$, we can select rectangles $R_-$ and $R_+$ satisfying \eqref{eq:eps_app} from the approximating set $\rr$ defined above.
\end{lemm}

To complete our characterization of the approximating set, it suffices to bound the cardinality of $\rr$. This computation is carried out in the Appendix, in the proof of Theorem \ref{theo:approx_set}.

 \subsection{Uniform Concentration over Rectangles}
\label{sec:concentration}

In the previous section, we showed how to $\varepsilon$-approximate all possible tree leaves under the Lebesgue measure on $[0, \, 1]^d$. However, to understand the behavior of decision trees, we do not want to approximate tree leaves in terms of Lebesgue measure, but rather in terms of the empirical measure induced by the training features $\{X_i\}_{i = 1}^n$, which, given our assumptions, are uniformly distributed over $[0, \, 1]^d$.

The following result lets us get around this issue by showing that the empirical measure induced by the training examples is concentrated enough that, with high probability, the set $\rrset$ is also a good approximating set in terms of the empirical measure induced by the training data. 

\begin{theo}
\label{theo:counts}

Suppose that Assumption \ref{assu:unif} holds, and that we have a sequence of problems indexed by $n$ with values of $d$ and $k$ satisfying Assumption \ref{assu:scaling}.
Let $\rrset$ be as defined in \eqref{eq:union} with $s$ as in \eqref{eq:s_def}, and choose $\varepsilon$ and $w$ such that
\begin{align}
\label{eq:wdef}
 \varepsilon = \frac{1}{\sqrt{k}}, \ \eqand \
w = \frac{1}{2 \zeta} \, \frac{k}{n},
\end{align}
where $\zeta \geq 1$ is the constant from Assumption \ref{assu:unif}.
Then, there exists an $n_0 \in \NN$ such that, for every $n \geq n_0$, the following statement holds with probability at least $1 - n^{-1/2}$:
for every possible leaf $L \in \law_{\alpha, \, k}$, we can select a rectangle $R \in \rrset$ such that $R \subseteq L$,
$\lambda\p{L} \leq e^{\varepsilon} \lambda\p{R}$,
and
\begin{align}
\label{eq:count_second}
\begin{split}
&\#L - \#R \leq 3 \, \zeta^2 \varepsilon \#L  + 2\sqrt{3 \log \p{\abs{\rrset}}\#L} + \oo\p{\log \p{\abs{\rrset}}}.
\end{split}
\end{align}
\end{theo}

We can then turn this result into a concentration bound on our empirical process of interest \eqref{eq:emp_proc},
provided we have a tail bound on the responses $Y_i$.
In the following result, we obtain such a tail bound by imposing a uniform sub-Gaussianity
requirement on the conditional distribution $Y \cond X \sim G\p{\cdot; \, X}$.
We define a conditional distribution $G\p{\cdot; X}$ to be uniformly sub-Gaussian
if there is a constant $M > 0$ for which the following holds.
For any $X$-marginal distribution $X \sim F_X(\cdot)$ supported on $[0, \, 1]^d$ and $t > 0$,
the resulting $Y$-marginal distribution $F_Y(\cdot) = \int G\p{\cdot; \, X} dF_X$ satisfies:
\begin{equation}
\label{eq:unif_subg}
\EE[Y \sim F_Y] {e^{t\p{Y - m}}} \leq e^{\frac{1}{2} \, t^2 M^2}, \ \where \ m = \EE[Y \sim F_Y]{Y}.
\end{equation}
Note that if $Y$ is bounded by $M$, i.e., $\abs{Y} \leq M$ almost surely, then \eqref{eq:unif_subg}
can immediately be verified using standard results.

\begin{coro}
\label{coro:rect}
Suppose that the conditions of Theorem \ref{theo:counts} hold, that the parameters $\varepsilon$ and $w$ are chosen as in \eqref{eq:wdef}, and that the conditional distribution of $Y$ given $X$ is uniformly sub-Gaussian in the sense of \eqref{eq:unif_subg}. Then, there exists an $n_0 \in \NN$ such that, for all $n \geq n_0$, the following holds with probability at least $1 - 2/\sqrt{n}$:
\begin{align}
\label{eq:rect1}
\begin{split}
& \sup_{L \in \law} \Bigg\{\abs{\frac{1}{\#L} \sum_{\cb{i : X_i \in L}} Y_i - \EE{Y | X \in L}}  \\
&\ \ \ \ \ \ \ \ - \p{2M + 4 \sqrt{3} \sup\cb{\abs{Y_i} : X_i \in L}} \, \sqrt{\frac{\log \p{\abs{\rrset}}}{k}} \Bigg\} \leq 0.
\end{split}
\end{align}
\end{coro}

\subsubsection{Proof Sketch for Theorem \ref{theo:counts}}
\label{sec:emp_conc}

Here, we present a series of technical results that lead up to Theorem \ref{theo:counts}. We begin with the following technical lemma, which follows from the Chernoff-Hoeffding concentration bound.
In practice, we will always use Lemma \ref{lemm:count_conc} with $\rr$ set to our finite approximating set $\rrset$.

\begin{lemm}
\label{lemm:count_conc}
Fix a sequence $\delta(n) > 0$, and define the event
\begin{align}
\label{eq:count_conc}
\mathcal{A} \ : \ &\sup \cb{\frac{\abs{\#R - n \, \mu(R)}}{\sqrt{n \, \mu\p{R}}} : R \in \rr, \, \mu\p{R} \geq \mu_{\min} }
\leq \sqrt{3 \log \p{\frac{\abs{\rr}}{\delta}}}
\end{align}
for any set of rectangles $\mathcal{R}$ and threshold $\mu_{\min}$.
Then, for any sequence of problems indexed by $n$ with
\begin{equation}
\label{eq:asymp_cond}
\limn \frac{\log \p{\abs{\rr}}}{n \, \mu_{\min}} = 0, \ \eqand \ \limn \frac{\delta^{-1}}{\abs{\rr}} = 0,
\end{equation}
there is a threshold $n_0$ such that, for all $n \geq n_0$, we have
$ \PP{\mathcal{A}} \geq 1 - \delta. $
Note that, above, $\mathcal{A}$, $\rr$, $\mu_{\min}$, and $\delta$ are all implicitly changing with $n$.
\end{lemm}

For now, the relation \eqref{eq:count_conc} is only valid for rectangles $R$ contained in $\rr$, which we will take to be our finite approximating set $\rrset$. In general, however, the leaves $L \in \law_{\alpha, \, k}$ we want to study will not be in $\rrset$. The following result lets us move beyond this issue by providing a bound that is valid for all rectangles $R$, not just those in $\rrset$.

\begin{lemm}
\label{lemm:all_count}
Suppose that the event $\mathcal{A}$ defined in Lemma \ref{lemm:count_conc} has occurred with $\rr = \rrset$ and $\mu_{\min} = \zeta w$, where $\zeta$ is the constant from Assumption \ref{assu:unif}. Then, all rectangles $R$ with $\mu\p{R} \geq \zeta w$ satisfy:
\begin{align*}
\#R \leq e^{\zeta^2 \varepsilon} n \mu\p{R} + e^{\frac{1}{2} \, \zeta^2 \varepsilon} \sqrt{3 n \mu\p{R} \log \p{\frac{\abs{\rrset}}{{\delta}}}}.
\end{align*}
\end{lemm}

Finally, in order for the above result to be useful for understanding the leaves of decision trees, we need to show that all possible leaves $L$ will satisfy the condition $\mu\p{L} \geq \zeta \, w$; the result below gives us such a guarantee.
With these results in hand, proving Theorem \ref{theo:counts} reduces to algebra.

\begin{coro}
\label{coro:size}
Let $\mathcal{A}$ be the event from Lemma \ref{lemm:count_conc} with $\rr = \rrset$ and $\delta = 1/\sqrt{n}$, and define the parameters $s$, $w$ and $\varepsilon$ as in the statement of Theorem \ref{theo:counts}. Then, there exists an $n_0 \in \NN$ such that, for $n \geq n_0$,
\begin{equation}
\label{eq:coro:size}
\inf_{L \in \law_{\alpha, \, k}} \cb{\mu\p{L}} \geq \zeta w \text{ on the event } \mathcal{A}.
\end{equation}
\end{coro}

\subsection{Proof of Theorem \ref{theo:main}}
\label{sec:main_result_proof}

We have now gathered all the ingredients required to prove Theorem \ref{theo:main}, which follows from combining Corollary \ref{coro:set_size_bound} with Corollary \ref{coro:rect}.
In fact, the first bound \eqref{eq:main_complete} follows directly from these two results by
noting that $\sup\cb{\abs{Y_i} \cond X_i \in L} \leq M$ (since $Y_i$ is bounded by
assumption) and that $2 + 4\sqrt{3} < 9$.

Now, if moreover $\liminf d/n > 0$, then we can verify that $\log(n)\log(d) - \log(n/k)\log(dk) \geq 0$ for large enough $n$, and so we can use Corollary \ref{coro:set_size_bound} to bound
\begin{equation}
\label{eq:rrset_bound}
 \log\p{\rr_{s, \, n, \, \varepsilon}} \leq \frac{\log\p{n} \log\p{d}}{\log\p{1 - \alpha}^{-1}} + \oo\p{\max\cb{\log \p{n} \log\p{\log\p{n}}, \, \log(d)}}.
\end{equation}
Thanks to our assumption that $\liminf d/n > 0$, the remained term is negligible and \eqref{eq:main} holds.

\section{Lower Bounds}
\label{sec:lower_bounds}

In this final section, we complement our main adaptive concentration bound, and show that the convergence rate given in \eqref{eq:main} cannot be improved. \citet{lin2006random} have also studied lower bounds for forest convergence; however, they only consider non-adaptive forests and so their lower bounds are substantially weaker. We show the following:

\begin{theo}
\label{theo:lowerbound}
For any $r >0$, set $d := d(n) = \lfloor n^r \rfloor$, and let $\alpha \leq 0.2$. Then, there exists a distribution over $(X, \, Y)$  and a sequence $k(n)$ satisfying the conditions Theorem \ref{theo:main} for which
\begin{align}
\label{eq:main2}
&\lim_{n \rightarrow \infty} \mathbb{P}\vast[\sup_{x \in [0, \, 1]^d, \, \Lambda \in \vv_{\alpha, \, k}} \left|T_\Lambda\p{x} - T^*_\Lambda\p{x} \right|
\geq \frac{M}{5} \,  \sqrt{\frac{ \log \p{n} \log \p{d}}{k}} \vast] = 1.
\end{align}
\end{theo}

Whenever, $\liminf d/n > 0$ (i.e., $r \geq 1$), the rate \eqref{eq:main2} has the same dependence on $n$, $d$, and $k$ as the upper bound \eqref{eq:main}, thus implying that our adaptive concentration bound from Theorem \ref{theo:main} is rate optimal.

To establish Theorem \ref{theo:lowerbound}, we take the $Y_i$ to be i.i.d. and independent of $\xx$, with 
$\PP{Y_1=M} = \PP{Y_1=-M} = 1/2$. We will construct $N=N(n)$ nodes $L_1,\ldots,L_N
\in \vv_{\alpha, \, k}\p{\xx}$ and then consider for $j=1,\ldots,N$:
\begin{equation}
\label{eq:tj_def}
T_j \ :=\ \frac{1}{\#L_j} \sum_{\{ i: X_i \in L_j\}} Y_i,
\ \eqand \
\tT_j \ :=\ \frac{1}{\#L_j} \sum_{\{ i: X_i \in L_j\}} \tY_i,
\end{equation}
where $\tT_j$ is an approximation to $T_t$ built using auxiliary random variables $\tY_i$ generated as
\begin{equation}
\label{eq:yi_approx}
\tY_i = Y_i \, \abs{Z_i} \where Z_i \simiid \nn\p{0, \, 1}.
\end{equation}
Notice in particular that the $\tY_i$ are jointly distributed
as independent Gaussian random variables with variance $M^2$.

The idea of the proof is to construct
a large set of candidate leaf nodes whose pairwise
intersections are small enough that the multivariate normal distribution of the
standardized $\tT_j$ has correlations that are bounded away from unity;
note that the leaves $L_j$ may be generated by different trees, and are thus allowed to overlap.
A normal approximation lemma \citep{leadbetter1983extremes} then allows us
to stochastically lower-bound the distribution of $\max_j \tT_j$
in terms of the distribution of a correlated multivariate normal that can be constructed in a simple way from
an i.i.d. normal sequence.
Specifically, we establish the following lower bound of the tail of the $\tT_j$.

\begin{lemm}
\label{lemm:GWA}
For any $r >0$, set $d := d(n) = \lfloor n^r \rfloor$, and let $\alpha \leq 0.2$. 
Then, there exists a sequence $k=k(n)$ satisfying Assumption \ref{assu:scaling}
and a set of $N$ $\{\alpha, \, k\}$-valid candidate leaf nodes $L_j$
chosen independently of the $Y_i$ and $\tY_i$,
for which
\begin{align}
\label{eq:01_lb}
&\lim_{n \rightarrow \infty} \PP{ \max_{j=1,\ldots,N} \tT_j\, \geq \, 
1.999 \, M \sqrt{\frac{2}{5}} \,\sqrt{\frac{ \log (N)}{k}} } \, =\, 1,  \eqand \\
&\log(N) =  \frac{\log \p{n} \log \p{d}}{\log(5)} \p{1 + o(1)}.
\end{align}
\end{lemm}

In the second step, we establish a coupling between $T_i$
and the $\tT_i$ that is tight enough to guarantee that the approximation error
$\max_j \{\tT_j -T_j\}$ is smaller than  $\max_j \, \tT_j$.
To get this coupling, we use the following bound on the moment-generating
function of $Y_i - \tY_i$.

\begin{lemm}
\label{lemm:GWLemma}
Let $\PP{Y=1}=\PP{Y=-1}=1/2$ and $Z \sim N(0,1)$ independent of $Y$. Then
\begin{equation}
\label{eq:GWLemma}
\EE{ \exp \cb{ t\p{Y-Y|Z|}}}\, \leq \, \exp\cb{\p{(1- \sqrt{\frac{2}{\pi}}} t^2}
\end{equation}
for $t$ in a neighborhood of zero.
\end{lemm}

This lemma implies the bound on  $\max_j \{\tT_j -T_j\}$ given below.
The lower bound in Theorem \ref{theo:lowerbound} then follows 
from Lemma \ref{lemm:GWA} and Corollary \ref{coro:GWB}, together with the observation
that $1.999 \, \sqrt{2/5} - 1 \geq \sqrt{\log 5}/5$.

\begin{coro}
\label{coro:GWB}
Suppose that the statistics $T_j$ and $\tT_j$ are constructed as in (\ref{eq:tj_def} - \ref{eq:yi_approx}), with leaves $L_j$ chosen independently of the $Y_i$ and $\tY_i$. Then,
\begin{equation}
\label{GWB}
\lim_{n \rightarrow \infty} \PP{ \max_{j=1, \, \ldots, \, N} \abs{\tT_j -T_j}\, \leq \, M \, \sqrt{\frac{\log\p{N}}{k}}} \, = \, 0.
\end{equation}
\end{coro}

\section*{Acknowledgment}

The authors are grateful for helpful feedback from two anonymous referees.

\bibliographystyle{plainnat}
\bibliography{references}

\printaddresses

\newpage

\begin{appendix}

\appendixtoc


\subsection*{Remark} In this appendix, we present our proofs in the order of logical dependence instead of the order in which they appear in the main text. For example, Theorem \ref{theo:approx_set} depends on Lemma \ref{lemm:eps_approx}, so we prove the lemma first.

\section{The Approximating Set of Rectangles}

\subsection{Proof of Lemma \ref{lemm:eps_approx}}

We focus on showing how to construct $R_+$; the construction of $R_-$ is analogous.
Recall that, given a rectangle
$$ R = \bigotimes_{j = 1}^s \left[r_j^-, \, r_j^+\right], 
\text{ our goal is to select a rectangle }
 R_+ = \bigotimes_{j = 1}^s \left[q_j^-, \, q_j^+\right] $$
from $\rr$ such that $R \subseteq R_+$ and
$ \lambda\p{R_+} \leq e^\varepsilon \lambda\p{R}.$
In order to guarantee this, it is sufficient to check that, for all $j$,
$$ q_j^- \leq r_j^-, \, r_j^+ \leq q_j^+, \eqand q_j^+ - q_j^- \leq e^{\varepsilon/s} \p{r_j^+ - r_j^-}. $$
Now, for each $j$, define
$$\tau_j = \left\lfloor \log_2 \frac{r_j^+ - r_j^-}{w} \right \rfloor, $$
let $q_j^-$ be the largest choice of the form \eqref{eq:counter} such that $q_j^- \leq r_j^-$, and pick $q_j^+$ analogously. These choices define a rectangle $R_+$ in $\rr$ such that $R \subseteq R_+$. By construction, we immediately see that
$$ 2^{\sum_{j = 1}^s \tau_j} \geq 2^{\sum_{j = 1}^s \p{\log_2 \frac{r_j^+ - r_j^-}{w} - 1}} \geq 2^{-s} w^{-\p{s - 1}}. $$
Moreover, by definition of $\tau_j$
$$ 2^{\tau_j} w \leq r_j^+ - r_j^- \leq 2^{\tau_j + 1} w, $$
thus, we can verify that
$$ q_j^+ - r_j^+, \, r_j^- - q_j^- \leq 2^{\tau_j - 1} \, \frac{w \varepsilon}{s} \leq \frac{1}{2} \frac{\varepsilon}{s} \p{r_j^+ - r_j^-}. $$
This implies that
$$ q_j^+ - q_j^- \leq \p{r_j^+ - r_j^-} \p{1 + \frac{\varepsilon}{s}}, $$
and so $|R_+| \leq e^\varepsilon |R|$.

\subsection{Proof of Theorem \ref{theo:approx_set}}

Given Lemma \ref{lemm:eps_approx}, in order to complete the proof of Theorem \ref{theo:approx_set} it suffices to bound the cardinality of the approximating set defined in Section \ref{sec:construction}. To do so, we first observe that for fixed values of $\{\tau_j\}$, the number of possible choices for the $\{a_j\}$ and $\{b_j\}$ is bounded by
\begin{align*}
\prod_{j = 1}^s &\p{1 + \left\lfloor 2^{1 - \tau_j} \frac{s}{w \varepsilon} \right\rfloor}\p{1 + \left\lceil \frac{2s}{\varepsilon} \right\rceil} \\
&= \p{\frac{4s^2}{w \varepsilon^2}}^s 2^{-\sum_{j = 1}^s \tau_j}  \cdot \p{1 + \oo\p{\varepsilon}}\\
&\leq \p{\frac{4s^2}{w \varepsilon^2}}^s 2^{s} \p{\prod_{j = 1}^s \frac{r_j^+ - r_j^-}{w}}^{-1}  \cdot \p{1 + \oo\p{\varepsilon}} \\
&= \frac{1}{w} \p{\frac{8s^2}{\varepsilon^2}}^s  \cdot \p{1 + \oo\p{\varepsilon}},
\end{align*}
because $\prod_{j = 1}^s \p{r_j^+ - r_j^-} \geq w$. Now, we can loosely bound the number of possible choices for $\{\tau_j\}$ by $\p{1 + \log_2 w^{-1}}^s$, yielding the desired bound.

\subsection{Proof of Corollary \ref{coro:set_size_bound}}
Given the parameter choices \eqref{eq:param_choices}, we can verify that
\begin{align*}
\log\abs{\rr_{s, \, w, \, \varepsilon}} & \leq \log \p{\binom{d}{s} \, \frac{1}{w} \p{\frac{8 s^2}{\varepsilon^2} \p{1 + \log_2 \left\lfloor \frac{1}{w} \right \rfloor}}^{s}  \cdot \p{1 + \oo\p{\varepsilon}}}\\
&=  \log \binom{d}{s} + 2 s \log\p{\varepsilon^{-1}} + 2 s \log(s)  \\
& \ \ \ \ \ \ \ \  + s\log\log\p{w^{-1}}+ \oo\p{\log\p{n}}.
\end{align*}
Meanwhile,
\begin{align*}
&\log \binom{d}{s} \leq s \log\p{d} = \frac{\log\p{n/k} \log\p{d}}{\log\p{\p{1 - \alpha}^{-1}}} + \oo\p{\log\p{d}}, \\
&s \log\p{\varepsilon^{-1}} =  \frac{1}{2} \, \frac{\log\p{n/k} \log\p{k}}{\log\p{\p{1 - \alpha}^{-1}}} + \oo\p{\log\p{n}}, \eqand\\
&s\log\p{s}, \ s\log\log\p{w^{-1}} = \frac{\log\p{n/k}\log\log(n)}{\log\p{\p{1 - \alpha}^{-1}}} + \oo\p{\log\p{n}}.
\end{align*}
Combining these results, we recover \eqref{eq:set_size_bound}.

\section{Concentration over Rectangles}

\subsection{Proof of Lemma \ref{lemm:count_conc}}

The proof of this result relies on a union bound. Our goal is to show that, for any rectangle $R$ with $R \in \rr$ and $\mu\p{R} \geq \mu_{\min}$, we can bound the large deviations of $\#R$ as follows: There is some $n_0 \in \NN$ such that, for all $n \geq n_0$,
\begin{align}
\label{eq:ub_target}
&\PP{\abs{\frac{\#R}{n} - \mu\p{R}} \geq \Delta} \leq \frac{\delta}{\abs{\rr}}, \text{ where}\\
&\Delta = \sqrt{\frac{\mu(R)}{n}} \, \sqrt{3 \log \p{\frac{\abs{\rr}}{\delta}}}.
\end{align}
Verifying \eqref{eq:ub_target} then immediately implies the desired bound on $\PP{\mathcal{A}}$.
We proceed in two parts. First, we verify that \eqref{eq:ub_target} holds for very large rectangles with $\mu(R) \geq 1/2$; second, we consider the smaller rectangles with $1/2 > \mu\p{R} \geq \mu_{\min}$.

In the case of large rectangles, we know that $\#R/n$ is sub-Gaussian with parameter
$\sigma_n^2 = 1/(4n)$. Thus,
\begin{align*}
&\PP{\abs{\frac{\#R}{n} - \mu\p{R}} \geq \Delta}
\leq 2 \exp\sqb{-2n \, \Delta^2} \\
&\ \ \ \ \ \ \leq 2 \exp\sqb{-3 \log \p{\frac{\abs{\rr}}{\delta}}} 
\leq 2 \p{\frac{\delta}{\abs{\rr}}}^3,
\end{align*}
and \eqref{eq:ub_target} is easily satisfied.

In order to analyze small rectangles, we need the tighter binomial concentration result of \citet{chernoff1952measure} and \citet{hoeffding1963probability}, stated below for convenience.

\begin{prop0}[Chernoff-Hoeffding]
Let $\#R$ be a binomial $(n, \, \mu)$ random variable. Then
\begin{align}
\label{eq:ch1}
&\PP{\frac{\#R}{n} \geq \mu + \Delta}  \leq \p{\p{\frac{\mu}{\mu + \Delta}}^{\mu + \Delta} \p{\frac{1 - \mu}{1 - \mu - \Delta}}^{1 - \mu - \Delta}}^n, \\
\label{eq:ch2}
&\PP{\frac{\#R}{n}  \leq \mu - \Delta}  \leq \p{\p{\frac{\mu}{\mu - \Delta}}^{\mu - \Delta} \p{\frac{1 - \mu}{1 - \mu + \Delta}}^{1 - \mu + \Delta}}^n.
\end{align}
\end{prop0}

Now by \eqref{eq:asymp_cond}, we know that $\Delta/\mu(R) \rightarrow 0$ uniformly over
our set of rectangles of interest.
Moreover, because we are working with small rectangles $R$, we also have $\Delta/(1 - \mu(R)) < \Delta/\mu(R)  \rightarrow 0$.
Finally, we can verify by calculus that
$$ \frac{1}{1 + x} \leq e^{-x + \frac{x^2}{2} + \abs{x}^3}\ \text{ for all } \ \abs{x} \leq 0.5. $$
Thus, for large enough $n$, the Chernoff-Hoeffding bound implies that for all
our rectangles of interest,
\begin{align*}
\PP{\frac{\#R}{n} \geq \mu + \Delta} &\leq \exp\bigg[ n \, \bigg( \p{-\frac{\Delta}{\mu} + \frac{\Delta^2}{2\mu^2} + \frac{\Delta^3}{\mu^3}}\p{\mu + \Delta} \\
&\hspace{-4em} + \p{\frac{\Delta}{1 - \mu} + \frac{\Delta^2}{2(1 - \mu)^2} + \frac{\Delta^3}{(1 - \mu)^3}}\p{1 - \mu - \Delta}\bigg)\bigg] \\
&\leq \exp\sqb{-   n \, \frac{\Delta^2}{2\mu(1 - \mu)} \p{1 + o(1)}}.
\end{align*}
We can also apply the same argument to \eqref{eq:ch2}, and get for large enough $n$:
\begin{align*}
\PP{\abs{\frac{\#R}{n} - \mu} \geq + \Delta}
&\leq 2  \exp\sqb{-   n \, \frac{\Delta^2}{2\mu(1 - \mu)} \p{1 + o(1)}} \\
&\leq \exp\sqb{-   n \, \frac{\Delta^2}{3\mu(1 - \mu)}} 
\leq \frac{\delta}{\abs{\rr}},
\end{align*}
thus concluding the proof.

\subsection{Proof of Lemma \ref{lemm:all_count}}

Because $\mu\p{R} \geq \mu_{\min}$, we can use Assumption \ref{assu:unif}
to verify that $\lambda\p{R} \geq \zeta^{-1} \mu\p{R} \geq w$.
Thus, by Theorem \ref{theo:approx_set}, we know that there exists a rectangles
$ R_+ \in \rrset $
such that
\begin{align*}
R \subseteq R_+, \eqand 
&e^{-\varepsilon} \lambda\p{R_+} \leq \lambda\p{R}.
\end{align*}
Moreover, again by Assumption \ref{assu:unif} and because $R \subseteq R_+$, we can verify that
\begin{align*}
\mu\p{R_+} &\leq \mu\p{R} + \zeta \p{\lambda\p{R_+} - \lambda\p{R}} \\
&\leq \mu\p{R} + \zeta\p{e^{\varepsilon} - 1} \, \lambda\p{R} \\
&\leq \mu\p{R} + \zeta^2\p{e^{\varepsilon} - 1} \, \mu\p{R} \\
&\leq e^{\zeta^2 \varepsilon} \mu\p{R}.
\end{align*}
Then, we see that on $\mathcal{A}$,
\begin{align*}
\#R
&\leq \#R_+ 
\leq n \mu\p{R_+} +  \sqrt{3 n \mu\p{R_+} \log \p{\frac{\abs{\rrset}}{\delta}}} \\
&\leq e^{\zeta^2  \varepsilon} n \mu\p{R} + e^{\frac{1}{2} \, \zeta^2 \varepsilon} \sqrt{3 n \mu\p{R} \log \p{\frac{\abs{\rrset}}{\delta}}},
\end{align*}
where the second inequality followed by Lemma \ref{lemm:count_conc}.

\subsection{Proof of Corollary \ref{coro:size}}

By Lemma \ref{lemm:all_count}, we see that on $\mathcal{A}$
\begin{align*}
&\sup\cb{\#R : \mu\p{R} = \zeta w}
\leq e^{\zeta^2 \varepsilon} \zeta n w + e^{\frac{1}{2} \, \zeta^2 \varepsilon} \sqrt{{3  \zeta n w \log\p{\frac{ \abs{\rrset}}{ \delta}}}} \\
&\ \ \ \ \ \ \ \ = e^{\zeta^2 \varepsilon} \frac{k}{2}   + e^{\frac{1}{2} \, \zeta^2 \varepsilon}\sqrt{\frac{3k}{2} \, \log\p{\frac{ \abs{\rrset}}{ \delta}}} < \frac{3k}{4}, 
\end{align*}
for large enough $n$.
In other words, for large enough $n$, all rectangles of size $\zeta w$ can have at most $3k/4$ points in them. Thus, we conclude that, on $\mathcal{A}$, any rectangle with $k$ points must have size greater than $\zeta w$.

\subsection{Proof of Theorem \ref{theo:counts}}

For this whole proof, we assume that the event $\mathcal{A}$ defined in Lemma
\ref{lemm:count_conc} has occurred, with $\rr = \rrset$, $\mu_{\min} = \zeta w/2$,
and $\delta = 1/\sqrt{n}$. Note that, thanks to Assumption \ref{assu:scaling},
these choices satisfy the conditions \eqref{eq:asymp_cond},
and so Lemma \ref{lemm:count_conc} implies that the is an $n_0 \in \NN$
such that the event $\mathcal{A}$ must occur
with probability at least $1 - 1/\sqrt{n}$ for all $n \geq n_0$.

Given the event $\mathcal{A}$, we first note that Corollary \ref{coro:size}
implies that
$\mu\p{L} \geq \zeta \, w$ for all $L \in \law_{\alpha, \, k}$;
by Assumption \ref{assu:unif}, this also implies that
$\lambda(L) \geq w$ for all $L \in \law_{\alpha, \, k}$.
Thus, by Theorem \ref{theo:approx_set}, for each possible leaf $L \in \law_{\alpha, \, k}$
we can select a rectangle $R \in \rrset$ such that $R \subseteq L$ and
$$ \lambda\p{L} \leq e^{\varepsilon} \lambda\p{R}. $$
This establishes the first part of our desired result.

Next, we need to control the counts $\#L$ and $\#R$ on $\mathcal{A}$.
First, by Lemma \ref{lemm:all_count}, we immediately see that,
on $\mathcal{A}$,
$$ e^{\zeta^2 \varepsilon} \, n \mu\p{L} + e^{\frac{1}{2} \, \zeta^2 \varepsilon}\sqrt{3\log \p{\frac{\abs{\rrset}}{{\delta}}}} \, \sqrt{ n \mu\p{L} } - \#L \geq 0, $$
or equivalently that
\begin{align*}
&n\mu\p{L} \geq \frac{e^{-\zeta^2 \varepsilon}}{4} \Bigg(-  \sqrt{3\log \p{\frac{\abs{\rrset}}{{\delta}}}} + \sqrt{3\log \p{\frac{\abs{\rrset}}{{\delta}}} + 4 \#L} \Bigg)^2,
\end{align*}
and so, because $\varepsilon \rightarrow 0$ and we are in a regime where
$\log \p{{\abs{\rrset}}/{{\delta}}} \ll \#L$, we find that
\begin{align*}
&e^{\zeta^2 \varepsilon} \, n\mu\p{L} \geq \#L - \sqrt{3\log \p{\frac{\abs{\rrset}}{{\delta}}}  \, \#L} \ + \oo\p{\log \p{\frac{\abs{\rrset}}{{\delta}}}}.
\end{align*}
Meanwhile, by Assumption \ref{assu:unif},
\begin{align*}
\mu(R) &\geq \mu(L) - \zeta \p{\lambda(L) - \lambda(R)}
\geq \mu\p{L} - \zeta \p{1 - e^{-\varepsilon}} \lambda(L) \\
&\geq \p{1 - \zeta^2 \p{1 - e^{-\varepsilon}}} \mu(L) \geq \mu_{\min},
\end{align*}
where the last inequality is valid
provided $\varepsilon$ is small enough (i.e., $n$ is large enough).
Thus, because $R \in \rrset$, we can use Lemma \ref{lemm:count_conc}
to verify that, on $\mathcal{A}$,
\begin{align*}
\#R
&\geq n\mu\p{R} - \sqrt{3 \, n\mu\p{R}\log \p{\frac{\abs{\rrset}}{{\delta}}}} \\
&\geq \p{1 - \zeta^2 \p{1 - e^{-\varepsilon}}} n \mu(L) - \sqrt{3n\mu\p{R}\log \p{\frac{\abs{\rrset}}{{\delta}}}} \\
&\geq e^{-2 \, \zeta^2 \varepsilon} n \mu(L) - \sqrt{3n\mu\p{R}\log \p{\frac{\abs{\rrset}}{{\delta}}}} \\
&\geq e^{-2 \, \zeta^2 \varepsilon} n \mu(L) - \sqrt{3\#R \log \p{\frac{\abs{\rrset}}{{\delta}}}} + \oo\p{\log \p{\frac{\abs{\rrset}}{{\delta}}}} \\
&\geq e^{-2 \, \zeta^2 \varepsilon} n \mu(L) - \sqrt{3\#L \log \p{\frac{\abs{\rrset}}{{\delta}}}} + \oo\p{\log \p{\frac{\abs{\rrset}}{{\delta}}}},
\end{align*}
where the third inequality relies on $\varepsilon$ being small enough, and the
fourth inequality relies on a second application of Lemma \ref{lemm:count_conc}.
Chaining these inequalities, together, we find that, on $\mathcal{A}$ and
provided that $n$ is large enough,
\begin{align*}
\#R
&\geq e^{-3 \, \zeta^2 \varepsilon} \#L - \p{1 + e^{-3 \, \zeta^2 \varepsilon}} \sqrt{3\#L \log \p{\frac{\abs{\rrset}}{{\delta}}}} \\
&\ \ \ \ \ \ \ + \oo\p{\log \p{\frac{\abs{\rrset}}{{\delta}}}}.
\end{align*}
Finally, by Assumption \ref{assu:unif},
$\varepsilon = 1/\sqrt{k} \ll 1/\sqrt{\log(|\mathcal{R}|/\delta)}$,
and moreover we know that $\delta^{-1} \ll \abs{\rrset}$;
thus, the expression simplifies to
\begin{align*}
&\#L - \#R \leq 3 \, \zeta^2 \varepsilon \#L + 2 \, \sqrt{3 \log \p{\frac{\abs{\rrset}}{{\delta}}} \, \#L}  + \oo\p{\log \p{\abs{\rrset}}}.
\end{align*}

\subsection{Proof of Corollary \ref{coro:rect}}

Let $\acal$ be the ``good'' set used in the proof of Theorem \ref{theo:counts}; recall that $\PP{\acal} \geq 1 - 1/\sqrt{n}$ for large enough $n$. 
Now, for any leaf $L$ generated by a valid tree, let $R \in \rrset$ be the inner approximation for $L$ constructed in Theorem \ref{theo:counts}. By the triangle inequality,
\begin{align*}
& \sup \cb{ \abs{\frac{1}{\#L} \sum_{\cb{i : X_i \in L}} Y_i - \EE{Y | X \in L}} : L \in \law} \\
& \ \ \ \ \leq \sup \Bigg\{\abs{\frac{1}{\#L} \sum_{\cb{i : X_i \in L}} Y_i  -  \frac{1}{\#R} \sum_{\cb{i : X_i \in R}} Y_i } : {L \in \law} \Bigg\} \\
& \ \ \ \  \ \ \ \ +  \sup \Bigg\{ \abs{\frac{1}{\#R} \sum_{\cb{i : X_i \in R}} Y_i - \EE{Y | X \in R}} :  R \in \rrset, \, \#R \geq k \Bigg\} \\
& \ \ \ \  \ \ \ \ + \sup \cb{\abs{\EE{Y | X \in R} - \EE{Y | X \in L}} : {L \in \law} }
\end{align*}
We can now proceed to bound each term individually.
Starting with the first term, we note that because $R \subseteq L$
\begin{align*}
&\abs{\frac{1}{\#L} \sum_{\cb{i : X_i \in L}} Y_i  -  \frac{1}{\#R} \sum_{\cb{i : X_i \in R}} Y_i }  \leq  2 \sup\cb{\abs{Y_i} : X_i \in L} \, \frac{\#L - \#R}{\#L},
\end{align*}
and by Theorem \ref{theo:counts}, on event $\acal$,
\begin{align*}
&\sup\cb{\frac{\#L - \#R}{\#L} : L \in \law} \leq 3 \, \zeta^2 \varepsilon + 2 \sqrt{\frac{3 \log\abs{\rrset}}{k}} + \oo\p{\frac{\log\abs{\rrset}}{k}}.
\end{align*}
Meanwhile, in order to bound the last term, we note that our uniform sub-Gaussianity condition
implies that, for any points $x_1$ and $x_2$, $|\EE{Y \cond X = x_2} - \EE{Y \cond X = x_1}| \leq 2M$.
To see this, apply \eqref{eq:unif_subg} with
$F(X) = \frac{1}{2} \delta\p{\cb{X = x_1}} + \frac{1}{2} \delta\p{\cb{X = x_2}}$. Thus, we find that
\begin{align*}
&\abs{\EE{Y | X \in R} - \EE{Y | X \in L}}  \leq  2M\,\frac{\mu\p{L} - \mu\p{R}}{\mu\p{L}} \leq 2M \zeta^2 \p{1 - e^{-\varepsilon}}, 
\end{align*}
where the last inequality used Assumption \ref{assu:unif} and the fact that $\lambda(R) \geq e^{-\varepsilon} \lambda(L)$.

Finally, thanks to our uniform sub-Gaussiantiy assumption \eqref{eq:unif_subg}, we see
that---conditionally on $\#R$---then mean of the $Y_i$ over $R$ is sub-Gaussian with
parameter $\sigma^2 = M^2/\#R$. Moreover, by the proof of Theorem \ref{theo:counts},
we see that the is an $n_0 \in \NN$ such that $\#R \geq k / 2$ on the $\mathcal{A}$
whenever $n \geq n_0$.
Thus, conditionally on $\mathcal{A}$ and provided that $n \geq n_0$, the following event $\bb$ holds with probability at least $1 - 1/\sqrt{n}$:
\begin{align*}
&\bb: \sup \Bigg\{ \abs{\frac{1}{\#R} \sum_{\cb{i : X_i \in R}} Y_i - \EE{Y | X \in R}}  \Bigg\}  \leq    M \sqrt{\frac{2 \, \log\p{2|\rr| \sqrt{n}}}{k/2}}.
\end{align*}
To verify this fact, we apply a union bound over the set of rectangles $R \in \rrset$
with $\#R \geq k / 2$.
Combining all these bounds yields \eqref{eq:rect1}.

\section{Random Forest Consistency}

\subsection{Proof of Lemma \ref{lemm:no_bad_split}}

By expanding the square, $4 \, N^- \, N^+ \,/\, (N^- + N^+)^2 \leq 1$, and so $\ell(\theta) \leq \Delta^2\p{\theta}$. Now, by Assumption \ref{assu:sparse}, for all $j \not\in\qq$,
$$ \Delta^*\p{\theta} := \EE{Y_i \cond X_i \in \nu, \, (X_i)_j > \theta} - \EE{Y_i \cond X_i \in \nu, \, (X_i)_j \leq \theta}  = 0. $$
Moreover, we constructed the tree such that the sets $\{i : X_i \in \nu, \, (X_i)_j \leq \theta\}$ and $\{i : X_i \in \nu, \, (X_i)_j > \theta\}$ both have at most $k$ observations. Thus, by Theorem \ref{theo:main} and Corollary \ref{coro:rect},
$$ \abs{\Delta\p{\theta} - \Delta^*\p{\theta}} = \abs{\Delta\p{\theta}} \leq 2 \times 9M \sqrt{\frac{\log\p{n}\log\p{d}}{k \log\p{\p{1 - \alpha}^{-1}}}}, $$
with probability at least $1 - \oo(1/\sqrt{n})$, uniformly over all possible nodes $\nu$ with at least $2k$ observations and all variables $j \not\in Q$. We conclude that, with probability tending to $1 - \oo(1/\sqrt{n})$, \eqref{eq:split_thresh} is never satisfied for any $j \not\in \qq$, uniformly over all nodes of all trees that can be generated as guess-and-check trees.

\subsection{Proof of Lemma \ref{lemm:good_split}}

Let $\nu$ be the current node considered by the algorithm, let $j \in \qq$, and suppose that the node $\nu$ has never yet been cut along the direction $j$. To establish the desired result, it suffices to show that the split will succeed with probability tending to 1, uniformly over all $j \in \qq$ and all possible nodes $\nu$ with at least $2k$ observations that have not yet been divided along $j$.
Our goal is to show that, with high probability, the split at $\theta = 1/2$ satisfies \eqref{eq:split_thresh}. Since the actual splitting point $\htheta$ is chosen by maximizing $\ell(\theta)$, the result then also holds for $\htheta$.

We first note that, by Assumption \ref{assu:monotone},
$$ \abs{\Delta^*\p{\frac{1}{2}}} = \abs{\EE{Y_i \cond X_i \in \nu, \, (X_i)_j > \frac{1}{2}} - \EE{Y_i \cond X_i \in \nu, \, (X_i)_j \leq \frac{1}{2}}} \geq \beta. $$
Thus, by Theorem \ref{theo:main} and Corollary \ref{coro:rect}, and by Assumption \ref{assu:scaling} on the minimum leaf-size $k$, we see that with probability at least $1 - \oo(1/\sqrt{n})$
$$ \Delta^2\p{\frac{1}{2}} \geq \beta^2 + o(1) $$
uniformly over all possible nodes.
Next, by Lemma \ref{lemm:all_count} and a similar argument, we again find that with probability at least $1 - \oo(1/\sqrt{n})$
$$ \frac{4 N_-\p{1/2} N_+\p{1/2}}{\p{N_-\p{1/2} + N_+\p{1/2}}^2} = 1 + o(1), $$
uniformly over all possible nodes. Thus, \eqref{eq:split_thresh} is in fact satisfied at $\theta =1/2$ with high probability, and the split will succeed with high probability (although the split may not actually occur at $\theta = 1/2$, as there may be even more significant potential splitting points).

\subsection{Proof of Theorem \ref{theo:unif_consistent}}

We start by defining the following event $\ee$: ``all trees in the forest never split on any variable $j \not\in \qq$, and always split on any variable $j \in \qq$ when $j$ is drawn in phase 2 of the guess-and-check procedure.'' On this event, our $d$-dimensional guess-and-check tree is equivalent to a $q$-dimensional guess-and-check tree supported on the coordinates $j \in \qq$.

From Corollary \ref{coro:main_rf}, we already know that
$$ \sup_{x \in [0, \, 1]^d} \abs{H(x) - H^*(x)} = \oo_p\p{\sqrt{\frac{\log(n)\log(d)}{k}}}, $$
where $H(x)$ is our guess-check-tree and $H^*(x)$ is the corresponding partition-optimal forest in the sense of definition \ref{defi:forest}. 
Moreover, from Lemma \ref{lemm:no_bad_split} and Lemma \ref{lemm:good_split}, we know that $\PP{\ee} \rightarrow 1$.
Thus, to obtain uniform consistency, it only remains to show that, conditionally on $\ee$,
$$ \sup_{x \in [0, \, 1]^d} \abs{H^*(x) - \EE{Y \cond X = x}} = o_p(1). $$
Now, let $T^*$ be a single partition-optimal tree comprising $H^*$. Because $\EE{Y \cond X = x}$ is Lipschitz in $x$, and because $T^*(x) = \EE{Y \cond X \in L(x)}$, we see that 
$$ \abs{T^*(x) - \EE{Y \cond X = x}} \leq C_{LIP} \, \diam\p{L(x)}, $$
where $C_{LIP}$ is the Lipschitz constant, and the diameter $\diam((L(x))$ of the leaf $L(x)$ is defined as the longest line segment contained inside $L(x)$. This implies that
$$ \abs{H^*(x) - \EE{Y \cond X = x}} \leq C_{LIP} \, \EE[*]{\diam\p{L(x)}}, $$
where $\mathbb{E}_*$ is an expectation over all trees comprising the forest. Finally, by Lemma 2 of \citet{meinshausen2006quantile}, we see that on event $\ee$
$$ \sup_{x \in [0, \, 1]^d} \EE[*]{\diam\p{L(x)}} = o_p(1), $$
thus concluding the proof. We note that, in his paper, \citet{meinshausen2006quantile} only discusses convergence at a single $x$; however, his proof is based on a Kolmogorov-Smirnov argument that in facts holds uniformly for all $x$.

\subsection{Proof of Theorem \ref{theo:l2_consistent}}

We again define the same event $\ee$ as in the proof of Theorem \ref{theo:unif_consistent}. This time, since the tree effectively always splits at the middle of a randomly selected feature $j \in \qq$, our tree on event $\ee$ is in fact equivalent to a {non-adaptive} median tree trained on the feature set $\qq$, where each splitting variable is selected independently at random. This is exactly the class of forests studied by \citet{duroux2016impact}, who showed that \eqref{eq:biau_rate} in fact holds for them.

Now, thanks to Lemma \ref{lemm:no_bad_split} and Lemma \ref{lemm:good_split}, the probability of the event $\ee$ failing decays to 0 at a rate $1/\sqrt{n}$. Thus, because our responses $Y$ are bounded, the risk accrued from failures of $\ee$ is bounded on the order of $\oo(1/\sqrt{n})$, which is vanishingly small relative to the rate in \eqref{eq:biau_rate}.
Thus, we conclude that \eqref{eq:biau_rate} in fact holds for guess-and-check trees with $\alpha = 1/2$.

\section{Lower Bounds}

\subsection{Proof of Lemma \ref{lemm:GWA}}

For any $0 < \alpha < 0.5$, we study $\alpha$-random
partitions generated as follows. Given $d_n = \lfloor n^r \rfloor$ as assumed in
Lemma \ref{lemm:GWA}, we set
$$
s=s_n:= \max\cb{0, \, \left\lfloor \log \p{\frac{\log^3 (n)}{n}}  (\log \alpha)^{-1} \right\rfloor}
\eqand
k=k_n:= \left\lfloor n \alpha^{s_n} \right\rfloor.
$$
Then $k \geq \log^3 (n)$ for all $n \geq 100$, so Assumption 2 is met.
Now, for each index set $S \subset \{1,\ldots,d\}$ with $|S|=s$, define the partition
$\Lambda_S$ as follows: For each $j = 1, \, ..., \, n$, in order, if $j \in S$,
then recursively split each leaf along the $j$-th feature in such a way that
each child-node contains at least a fraction
$\alpha$ of the data points in its parent-node, with one child-node containing exactly
a fraction $\alpha$ up to rounding.
Given this construction, each terminal node $L$ has at least
$\#L \geq \alpha^s n \geq k $, so
$\Lambda_S$ is in fact an $\{\alpha, k\}$-valid partition.

Furthermore, the above construction provides for one terminal node $L$ that satisfies
\begin{align*}
\#L \leq \alpha^s n + \sum_{i=0}^{s-1} \alpha^i \leq  \p{k + 1} + \frac{1}{1-\alpha}  \leq k + 3.
\end{align*}
For each $s$-combination $S$ of $\{1,\ldots,d\}$ construct $\Lambda_S$ as described
above. Then for each of the $N:=\binom{d}{s}$ $s$-combinations the resulting
partition has one terminal node with the above properties. Denote these $N$ terminal
nodes by $L_1, \, \ldots, \, L_N$; so for $i=1, \, \ldots, \, N$, $k \leq \#L_i \leq k+3$.
Moreover, if $i \neq j$ then the splits in $L_i$ and $L_j$ occur on axes that differ
in at least one index. Since $X$ has independent marginals we get
\begin{equation}
\label{GW2}
\EE{ \#(L_i \cap L_j)}\ \leq \alpha \EE{\#L_i} +1 \ \leq \ \alpha k+3.
\end{equation}
Since the overlap between the $L_j$ is not very large, we might hope that the
maximum of the $\tT_j$ would be of comparable size to the
maximum of $N$ independent Gaussian random variables with variance
$\#L_j^{-1} M^2$. The following sub-result shows that this intuition is valid,
at least to within a factor $\sqrt{1 - \alpha}$.

\begin{prop0}
Assume that $d_n = n^r$ for some $r > 0$.
Then, given the nodes $L_j$ and statistics $\tT_j$ as constructed above, and for any $\eta > 0$,
\begin{equation}
\label{GWA}
\lim_{n \rightarrow \infty} \PP{ \max_{j=1,\ldots,N} \tT_j\ \geq \ 
(1 - \eta) \, M \, \sqrt{2(1-\alpha)} \sqrt{\frac{ \log (n) \log (d)}{\log \bigl(\alpha^{-1}\bigr)}}
\frac{1}{\sqrt{k}} } \ =\ 1.
\end{equation}
\end{prop0}

The claim from Lemma \ref{lemm:GWA} follows directly from \eqref{GWA} by noting that,
following our construction, an $\tilde{\alpha}$-random partition is always
an $\cb{\alpha, \, k}$-valid partition for any $\tilde{\alpha} \geq \alpha$. Thus,
if we know that $\alpha \leq 0.2$, we can apply \eqref{GWA} for 0.2-random partitions,
thus yielding the desired bound.

We now proceed to proving the sub-result.
Independence of $X$ and $\tY$ implies that
the $\sqrt{\#L_j}\ \tT_j$ are standard normal, and so, for $j \neq m$:
\begin{align*}
 \Cov{\sqrt{\#L_j} \ \tT_j, \, \sqrt{\#L_m} \ \tT_m } 
&= M \EE{\frac{\#(L_j \cap L_m)}{\sqrt{\#L_j \#L_m}}} \\
&\hspace{-3em}\leq M \EE{\frac{\#(L_j \cap L_m)}{k}} \leq M \p{\alpha + \frac{3}{k}}
\end{align*}
by \eqref{GW2}.
Now, let $Z_1, \, \ldots, \, Z_{N+1}$ be i.i.d. $\nn(0, \, M^2)$ and for $j=1, \, \ldots, \, N$ set
$$
\tZ_j \ :=\ \sqrt{1-\alpha_n}\, Z_j +\sqrt{\alpha_n}\,Z_{N+1}, \ \text{ with } \ \alpha_n := \alpha + \frac{3}{k_n}.
$$
The $\tZ_j$ are marginally normal with variance $M^2$, and $\text{Cov}[\tZ_j, \, \tZ_m]=\alpha_n M^2$
if $j \neq m$. Using Corollary~4.2.3
to the normal approximation lemma of \citet{leadbetter1983extremes},
and the fact that $\sqrt{\#L_j/k} \leq 1 + 2/k$,
we get for every $u>0$:
\begin{align}
\label{eq:C1}
\begin{split}
\PP{\max_{j=1,\ldots,N} \sqrt{k} \, \tT_j  \leq u} & \leq \PP{
  \max_{j=1,\ldots,N} \sqrt{\#L_j}\,\tT_j \leq \p{1+\frac{2}{k}}u} \\
& \leq \PP{\max_{j=1,\ldots,N} \tZ_j \leq \p{1+\frac{2}{k}}u} \\
& = \PP{\max_{j=1,\ldots,N} Z_j \leq \frac{(1+\frac{2}{k})u -\sqrt{\alpha_n} Z_{N+1}}{
  \sqrt{1-\alpha_n}}}.
  \end{split}
\end{align}
Setting $u= u_n:=(1-\eta)M\sqrt{2(1-\alpha)} \sqrt{{\log (n) \log (d)}/{\log(\alpha^{-1})}}$, our goal is to show that the above probability converges to 0.
First, observe that
\begin{equation*}
\PP{\frac{2u_n}{k_n} -\sqrt{\alpha_n} Z_{N+1} \leq \frac{\eta/2}{1 - \eta} u_n } \rightarrow 1
\end{equation*}
because, by assumption, $d_n = \lfloor n^r \rfloor$ whereas $k_n \sim \log(n)^3$.
Thus, by \eqref{eq:C1}, it suffices to verify that
\begin{equation}
\label{eq:Cgoal}
\PP{\max_{j=1,\ldots,N} Z_j \leq  \frac{1 - \eta/2}{1 - \eta} \frac{u_n}{\sqrt{1 - \alpha_n}}} \rightarrow 0.
\end{equation}
Then, noting that
$$\log(N) = \log\binom{d}{s} = s \log(d) \p{1 + o(1)} = \frac{\log(n)\log(d)}{\log\p{\alpha^{-1}}}\p{1 + o(1)}, $$
we can use a standard Gaussian tail bound to check that \eqref{eq:Cgoal} in fact holds.

\subsection{Proof of Lemma \ref{lemm:GWLemma}}
One readily checks that 
$$
\EE{ \exp \Bigl\{ t\bigl(Y-Y|Z|\bigr)\Bigr\}} = \exp \Bigl(\frac{1}{2} t^2 +t\Bigr)
\Phi(-t) +\exp \Bigl(\frac{1}{2} t^2 -t\Bigr) \Phi(t),
$$
where $\Phi$ is the cdf of $Z$. Using the expansion 
$$\Phi(t) =\frac{1}{2}+ \frac{1}{\sqrt{2\pi}} \sum_{k=0}^{\infty} \frac{(-1)^k t^{2k+1}}{
k! 2^k (2k+1)}, $$
we find that
\begin{align*}
& e^t \Phi(-t) + e^{-t} \Phi(t) \\
& = \cosh (t) -\sqrt{\frac{2}{\pi}} \, \sinh (t) \, 
\sum_{k=0}^{\infty} \frac{(-1)^k t^{2k+1}}{k! 2^k (2k+1)} \\
& = \sum_{k=0}^{\infty} \frac{t^{2k}}{(2k)!} \ -\ \sqrt{\frac{2}{\pi}}
  \p{ \sum_{k=0}^{\infty} \frac{t^{2k+1}}{(2k+1)!} } \p{
  \sum_{k=0}^{\infty} \frac{(-1)^k t^{2k+1}}{k! 2^k (2k+1)} } \\
& = 1\ +\ \p{\frac{1}{2} - \sqrt{\frac{2}{\pi}}}\, t^2 \ +\ \frac{1}{4\,!}\, t^4
  \ +\ O(t^6).
\end{align*}
Lemma~\ref{lemm:GWLemma} follows since
$$\frac{1}{4\,!} < {\p{\frac{1}{2} - \sqrt{\frac{2}{\pi}}}^2} \, \big/ \, {2\,!}. $$

\subsection{Proof of Corollary \ref{coro:GWB}}
For simplicity, we take $M = 1$, and so $Y_i \in \pm 1$ and $\text{Var}[\tY_i] = 1$.
Q standard argument using Markov's inequality gives for any $t,v>0$:
\begin{align*}
&\PP{\sqrt{\#L_j} \p{\tT_j-T_j}>v \Big|X} \\ 
&\ \ \ \ \ \ \leq \p{ \EE{\exp \cb{ 
t \p{Y_1 |Z_1| -Y_1 } }}}^{\#L_j} \exp \cb{- \sqrt{\#L_j}\,tv}\\
&\ \ \ \ \ \  \leq \exp \cb{ - \,\frac{v^2}{4 \p{1-\sqrt{\frac{2}{\pi}}}} }
\end{align*}
by Lemma~\ref{lemm:GWLemma}, provided that $t:={v}/\bigl(2\sqrt{\#L_j}\bigl(1-
\sqrt{{2}/{\pi}}\bigr)\bigr)$ is small enough. Now, set $v=v_n:=\eta \sqrt{\log N}$
for some fixed $\eta >0$. First, recalling $\min_j \#L_j \geq k_n \geq \log^3 n$,
we verify that in fact $t \rightarrow 0$. Thus,
\begin{align*}
& \PP{ \max_{j=1,\ldots,N} \sqrt{k_n} \p{\tT_j-T_j} > \eta \sqrt{ \log N}}\\
&\ \ \ \ \leq N \max_j \EE{ \PP{ \sqrt{\#L_j} \p{\tT_j-T_j} >v_n \Big|X} } \\
&\ \ \ \ \leq \exp \cb{ \bigl( \log N \bigr) \p{1-\frac{\eta^2}{4
  \p{1-\sqrt{{2}/{\pi}}}} }},
\end{align*}
which converges to zero provided that $\eta^2 > 4\bigl(1-\sqrt{{2}/{\pi}}\bigr) 
> 0.8$.

\end{appendix}

\end{document}